\documentclass[11pt]{article}
\newcommand{\documentdate}{6 VI 2023}

\usepackage{a4wide,latexsym,amsmath,varioref,graphicx,xcolor,amssymb,diagbox}

\title{Complexity of a Class of First-Order Objective-Function-Free Optimization
  Algorithms}

\author{
   S. Gratton%
   \thanks{Universit\'e de Toulouse, INP, IRIT, Toulouse, France. Email:
     serge.gratton@enseeiht.fr. Work partially supported by 3IA Artificial and
     Natural Intelligence Toulouse Institute (ANITI), French "Investing for the Future
     - PIA3" program under the Grant agreement ANR-19-PI3A-0004"}, 
   ~S. Jerad%
   \thanks{ANITI, Universit\'e de Toulouse, INP, IRIT, Toulouse, France. Email:
     sadok.jerad@enseeiht.fr}
   ~and Ph. L. Toint%
   \thanks{NAXYS, University of Namur, Namur, Belgium. Email:
     philippe.toint@unamur.be.
     Partially supported by ANITI.}
}

\newcommand{\beqn}[1]{\begin{equation}\label{#1}}
\newcommand{\eeqn}{\end{equation}}
\newcommand{\req}[1]{(\ref{#1})}
\newcommand{\ms}{\;\;\;\;}
\newcommand{\tim}[1]{\;\; \mbox{#1} \;\;}

\newtheorem{theorem}{Theorem}[section]
\newtheorem{lemma}[theorem]{Lemma}

\newcommand{\numsection}[1]{\section{#1}\setcounter{equation}{0}}
\newtheorem{corollary}[theorem]{Corollary}

\newcounter{algo}[section]
\renewcommand{\thealgo}{\thesection.\arabic{algo}}
\newcommand{\llem}[2]{\vspace{\baselineskip} 
\noindent\framebox[\textwidth]{\parbox{0.95\textwidth}{
\begin{lemma} \label{#1} \rm #2 \end{lemma} } } \vspace{\baselineskip} }

\newcommand{\algo}[3]{\refstepcounter{algo}
\begin{center}\begin{figure}[htbp]
\framebox[\textwidth]{
\parbox{0.95\textwidth} {\vspace{\topsep}
{\bf Algorithm \thealgo : #2}\label{#1}\\
\vspace*{-\topsep} \mbox{ }\\
{#3} \vspace{\topsep} }}
\end{figure}\end{center}}
\newcommand{\bpr}{{\bf Proof.} \hspace{1.5mm}}
\newcommand{\epr}{\hfill $\Box$ \vspace*{1em}}
\newcommand{\proof}[1]{
\begin{list}{}{
\setlength{\topsep}{0.0pt}
\setlength{\partopsep}{0.0pt}
\setlength{\leftmargin}{0.025\textwidth}
\setlength{\rightmargin}{0.5\leftmargin}
\setlength{\labelwidth}{0.5\leftmargin}
\setlength{\labelsep}{0.25\leftmargin}}
\item \bpr #1 \epr \noindent
\end{list}}
\newcommand{\lthm}[2]{\vspace{\baselineskip} 
\noindent\framebox[\textwidth]{\parbox{0.95\textwidth}{
\begin{theorem} \label{#1} \rm #2 \end{theorem} } } \vspace{\baselineskip} }
\newcommand{\ii}[1]{\{ 1, \ldots, #1 \}}
\newcommand{\iiz}[1]{\{ 0, \ldots, #1 \}}
\newcommand{\iibe}[2]{\{ #1, \ldots, #2 \}}

\newcommand{\calO}{{\cal O}}

\renewcommand{\Re}{\hbox{I\hskip -2pt R}}
\newcommand{\smallRe}{\hbox{\footnotesize I\hskip -2pt R}}
\newcommand{\bigfrac}[2]{\frac{\displaystyle #1}{\displaystyle #2}}
\newcommand{\bigsum}{\displaystyle \sum}

\newcommand{\bigmax}{\displaystyle \max}
\newcommand{\sfrac}[2]{{\scriptstyle \frac{#1}{#2}}}
\newcommand{\kap}[1]{\kappa_{\mbox{\tiny #1}}}
\newcommand{\eqdef}{\stackrel{\rm def}{=}}
\newcommand{\al}[1]{{\footnotesize{\sf #1}}}

\newcommand{\tal}[1]{{\normalsize {\sf #1}}}
\newcommand{\half}{\sfrac{1}{2}}

\newcommand{\flow}{f_{\rm low}}
\newcommand{\sgn}{{\rm sign}}
\newcommand{\neol}{\nonumber\\}
\newcommand{\kB}{\kap{B}}
\newcommand{\kBBL}{\kap{BBL}}
\DeclareMathOperator*{\average}{average}
\newcommand{\comment}[1]{}

\date{\documentdate}

\begin{document}

\maketitle

\begin{abstract}
A parametric class of trust-region algorithms for unconstrained nonconvex 
optimization is considered where the value of the objective function is never 
computed. The class contains a deterministic version of the first-order 
Adagrad method typically used for minimization of noisy function, but also 
allows the use of (possibly approximate) second-order information when available.
The rate of convergence of methods in the class is analyzed and is shown to be
identical to that known for first-order optimization methods using both
function and gradients values, recovering existing results for
purely-first order variants and improving the explicit dependence on problem 
dimension. This rate is shown to be essentially sharp. 
A new class of methods is also presented, for which a slightly worse and 
essentially sharp complexity result holds. Limited numerical experiments 
show that the new methods' performance may be comparable to that of standard 
steepest descent, despite using significantly less information, and that this 
performance is relatively insensitive to noise.
\end{abstract}


\numsection{Introduction}

This paper is concerned with
\underline{O}bjective-\underline{F}unction-\underline{F}ree
\underline{O}ptimization (OFFO) algorithms, which we define as numerical
optimization methods in which \emph{the value of the problem's objective 
function is never calculated}, although we obviously assume that it
exists.  This is clearly at variance with the large majority of available
numerical optimization algorithms, where the objective function is typically
evaluated at every iteration, and its value then used to assess progress
towards a minimizer and (often) to enforce descent. Dispensing with this
information is therefore challenging. As it turns out, first-order OFFO
methods (i.e. OFFO methods using gradients only) already exist for some time
and have proved popular and useful in fields such as machine learning or
sparse optimization, thereby justifying our interest. Many algorithms have 
been proposed, regrouped under the denomination of adaptive gradient methods, such as
Adagrad \cite{McMaStre10,DuchHazaSing11},
RMSprop \cite{TielHint12},
ADADELTA \cite{Zeil12},
Adam \cite{KingBa15},
SC-Adagrad \cite{MukkHein17},
WNGrad \cite{WuWardBott18}
or AMSgrad \cite{ReddKaleKuma18}
to cite a few. 

Adagrad was originally proposed for online learning \cite{McMaStre10,DuchHazaSing11}, 
where strong convergence guarantee were derived using current and past gradient 
information for a class of problems where the objective function changes at every iteration. 
The component-wise weighting of the gradient descent direction described in 
\cite{DuchHazaSing11} has proved to be remarkably successful in practice and has prompted 
the proposal of many variants. Among those, RMSprop \cite{TielHint12}, ADADELTA \cite{Zeil12}, 
Adam \cite{KingBa15} use an exponentially decreasing moving average when updating 
the weights instead of the non-decreasing technique proposed in \cite{DuchHazaSing11}. 
In particular, Adam \cite{DuchHazaSing11} has exhibited excellent performance 
in deep learning applications, as illustrated in recent numerical survey 
\cite{SchmSchnHenn21}. Closely related variants were also shown to be practically competitive 
in  \cite{ChenGu18}. However, its potential divergence was shown in \cite{ReddKaleKuma18}
where AMSGrad, a non-decreasing convergent scheme, was proposed as an alternative. The 
original theoretical analysis of Adagrad \cite{DuchHazaSing11,McMaStre10} was subsequently
refined for strongly convex functions \cite{MukkHein17} with a slight modification of 
the update rule, extended for minmax problems \cite{YangLiHe22}, further adapted 
for online use \cite{OrabPal15}, modified to ensure privacy \cite{Asietal21},
enhanced with a linesearch  for specific machine learning problems \cite{Vaswanietal20}, 
or with an accelerated-gradient-like update \cite{LevyAlpVolk18}. Of special interest is 
WNGrad \cite{WuWardBott18} which uses a mechanism to adapt the learning rate (or stepsize) 
using gradient information. Their update rule is very close to Adagrad-Norm 
\cite{WardWuBott19}, where a global (rather than componentwise) weight is used for all 
components. 
However, and despite its excellent performance, an analysis of Adagrad in the 
nonconvex setting appeared only relatively recently in \cite{LiOrab19} and \cite{WardWuBott19}. 
\cite{LiOrab19} requires the \textit{a-priori} knowledge of the Lipschitz constant (a serious theoretical drawback)
whereas \cite{WardWuBott19} is parameter agnostic. Both assume boundedness 
of the gradients for stochastic problems and obtain the conclusion that Adagrad's 
global rate of convergence is comparable to that of well-tuned stochastic gradient methods 
\cite{GhadLan13} up to logarithmic factors.  Other proofs, still requiring a 
uniform bound on the sampled gradient, were subsequently developed using simpler arguments 
\cite{DefoBottBachUsun22}, yielding improved convergence rates under 'gradient sparsity' 
\cite{ZhouChenTangYangCaoGu20}. The assumption of bounded stochastic gradients was 
finally removed, see \cite{Fawetal22,FawRoutCaraShak23,AttiKore23} and the references therein. 
A complexity analysis of Adagrad-Norm for deterministic nonconvex optimization was 
established in \cite{WardWuBott19}, showing a global rate of convergence equal to that 
of standard exact first order methods, without assuming bounded gradients. 
Inspired by this reference, \cite{TraoPauw21} 
proposed a new analysis of a (componentwise) Adagrad for the convex case, which features
an explicit\footnote{That is ignoring the potential dependence of the problem Lipschitz 
constant on dimension.} dependence on the problem's size which does not appear in the 
analysis of \cite{WardWuBott19}. 

In this paper, we consider the deterministic Adagrad algorithm \cite{DuchHazaSing11} 
using a trust-region point \cite{ConnGoulToin00} of view. 
Although this approach has not been used much in the machine learning context, 	
trust region techniques have been investigated in the framework of noisy optimization problems.
The method proposed in \cite{FanYuan01} uses function values but updates the radius 
bounding the steplength as a function of the gradient norm, which is also the case 
of \cite{CurtLubbRobi18}. The algorithm of \cite{BlanCartMeniSche19} 
and \cite{BellGuriMoriToin21b} also uses the (noisy) objective function's value 
while also allowing noise in the derivatives. In contrast, the algorithms described in \cite{GrapStel22} 
and \cite{CurtScheShi19} do not evaluate the objective function. They differ from the 
proposal we are about to describe both in the technique of proof and the fact that they do 
not subsume Adagrad or many of its variants. Moreover, the analysis of \cite{CurtScheShi19} requires the explicit 
knowledge of the problem's Lipschitz constant in the algorithm for obtaining the best complexity estimate.

The purpose of the present paper is to bridge the gap between standard first-order OFFO 
methods such as Adagrad and OFFO trust-region algorithms by considering a unified 
framework.  More specifically,
\begin{enumerate}
\item we re-interpret the deterministic Adagrad as a particular member of a fairly general 
parametric class of trust-region methods (Sections~\ref{method} and \ref{adap-s}). This 
class not only contains purely first-order algorithms such as Adagrad, 
but also allows the use of (possibly approximate) second-order information, should it be 
available using a Barzilai-Borwein approach \cite{BarzBorw88}, a limited-memory BFGS 
technique \cite{LiuNoce89} or even exact second derivatives.
\item We then provide, for our proposed class, an essentially sharp global\footnote{I.e., valid at every 
iteration.} bound on the gradient's norm as a function of the iteration counter, 
which is identical to that known for first-order optimization methods using both 
function and gradients values. This complexity result does not assume bounded gradients and 
extends that of \cite{WardWuBott19} only valid for Adagrad-Norm to the complete class.  
It also uses one of the available parameters of the class to mitigate the explicit 
dependence of that bound on the problem's dimension.
\item We next exploit the proposed OFFO trust-region framework of Section~\ref{method} 
to propose (in Section~\ref{newclass-s}) a new class class of 
such methods, for which an essentially sharp complexity result is also provided.
\item We finally illustrate our proposals by discussing 
some numerical experiments in Section~\ref{numerics-s}, suggesting that the considered OFFO methods may indeed be 
competitive with steepest descent in efficiency and reliability while being 
much less sensitive to noise.
\end{enumerate}

\noindent
\textbf{Notations.} In what follows, the superscript $T$ denotes the
transpose and $w_{i,k}$ denotes the $i$-th component of a vector
$w_k\in \Re^n$. Unless specified otherwise, $\|\cdot\|$ is the
Euclidean norm on $\Re^n$.  We also say, for non-negative quantities
$\alpha$ and $\beta$, that $\alpha$ is $\calO(\beta)$ is there exists
a finite constant $\kappa$ such that $\alpha \leq \kappa \beta$. 

\numsection{A class of first-order minimization methods}
\label{method}

We consider the problem
\beqn{problem}
\min_{x\in\smallRe^n} f(x)
\eeqn
where $f$ is a smooth function from $\Re^n$ to $\Re$. In particular, we will
assume in what follows that
\vspace*{2mm}
\noindent
\begin{description}
\item[AS.1:] the objective function $f(x)$ is continuously differentiable;
\item[AS.2:] its gradient $g(x) \eqdef \nabla_x^1f(x)$ is Lipschitz continuous with
   Lipschitz constant $L\geq 0$, that is
   \[
   \|g(x)-g(y)\| \le L \|x-y\|
   \]
   for all $x,y\in \Re^n$; 
\item[AS.3:] there exists a constant $\flow$ such that, for all $x$, $f(x)\ge \flow$.
\end{description}

\noindent
AS.1, AS.2 and AS.3 are standard for the complexity analysis of optimization
methods seeking first-order critical points, AS.3 guaranteeing in particular
that the problem is well-posed. We stress once more that we do not
assume that the gradient are uniformly bounded, at variance with
\cite{WuWardBott18,DefoBottBachUsun22,ZhouChenTangYangCaoGu20,GratJeraToin22a}. 

The class of methods of interest here are iterative and generate a sequence of
iterates $\{x_k\}_{k\geq0}$. The move from an iterate to the next directly
depends on the gradient at $x_k$ and algorithm-dependent \emph{scaling factors}
$\{w_k = w(x_0, \ldots, x_k)\}$ whose main purpose is to control the move's magnitude.
In our analysis, we will assume that
\begin{description}
  \item[AS.4:] for each $i\in\ii{n}$ there exists a constant $\varsigma_i\in(0,1]$
    such that, $w_{i,k} \geq \varsigma_i$ for all $k\geq0$,
\end{description}

Since scaling factors are designed to control the length of the step, they are
strongly reminiscent of the standard mechanism of the much studied trust-region optimization
methods (see \cite{ConnGoulToin00} for an extensive coverage and \cite{Yuan15}
for a more recent survey). In trust-region algorithms, a model of the objective
function at an iterate $x_k$ is built, typically using a truncated Taylor
series, and a step $s_k$ is chosen that minimizes this model with a
\emph{trust-region}, that is a region where the model is assumed to represent
the true objective function sufficiently well. This region is a ball around
the current iterate, whose radius is updated adaptively from iteration to
iteration, based on the quality of the prediction of the objective function
value at the trial point $x_k+s_k$. For methods using gradient only, the model
is then chosen as the first two terms of the Taylor's expansion of $f$ at the
iterate $x_k$. Although, we are interested here in methods where the objective
function's value is not evaluated, and therefore cannot be used to accept/reject
iterates and update the trust-region radius, a similar mechanism may be
designed, this time involving the weights $\{w_k\}$, the choice of which will 
detailed in the following two sections for two algorithmic subclasses of interest.
The resulting algorithm, which we call \al{ASTR1} (for Adaptively Scaled Trust Region using
1rst order information)  is stated \vpageref{ASTR1}.

\algo{ASTR1}{\tal{ASTR1}}
{
\begin{description}
\item[Step 0: Initialization. ]
  A starting point $x_0$ is given. Constants $\kB  \ge 1$ and $\tau \in (0,1]$  are also
  given.   Set $k=0$.
\item[Step 1: Define the TR. ]
  Compute $g_k = g(x_k)$ and define
  \beqn{Delta-def-a}
  \Delta_{i,k}  =\frac{|g_{i,k}|}{w_{i,k}}
  \eeqn
  where $w_k=w(x_0,\ldots,x_k)$.
\item[Step 2: Hessian approximation. ]
  Select a symmetric Hessian approximation $B_k$ such that
  \beqn{Bbound-a}
  \|B_k\| \le \kB. 
  \eeqn
\item[Step 3: GCP.]
  Compute a step $s_k$ such that
  \beqn{sbound-a}
  |s_{i,k}| \le \Delta_{i,k}  \ms (i \in \ii{n}),
  \eeqn
  and
  \beqn{GCPcond-a}
  g_k^Ts_k + \half s_k^TB_ks_k \le \tau\left(g_k^Ts_k^Q + \half (s_k^Q)^TB_ks_k^Q\right),
  \eeqn
  where
  \beqn{sL-def-a}
  s_{i,k}^L = -\sgn(g_{i,k})\Delta_{i,k},
  \eeqn
  \beqn{sQ-def-a}
  s^Q_k = \gamma_{k} s_k^L,
  \eeqn
  with
  \beqn{gammak-def-a}
  \gamma_k =
  \left\{ \begin{array}{ll}
  \min\left[ 1, \bigfrac{|g_k^Ts_k^L|}{(s_k^L)^T B_k s_k^L}\right] & \tim{if }
  (s_k^L)^T B_k s_k^L > 0,\\
  1 & \tim{otherwise.}
  \end{array}\right.
  \eeqn

\item[Step 4: New iterate.]
  Define
  \beqn{xupdate-a}
  x_{k+1} = x_k + s_k,
  \eeqn
  increment $k$ by one and return to Step~1.
\end{description}
}

\noindent
The algorithm description calls for some comments.
\begin{enumerate}
\item Observe that we allow the use of second-order information by
  effectively defining a quadratic model
  \beqn{q-model}
  g_k^Ts + \half s^TB_ks
  \eeqn
  where $B_k$ can of course be chosen as the true second-derivative matrix of
  $f$ at $x_k$ (provided it remains bounded to satisfy \req{Bbound-a}) or any
  approximation thereof.  Choosing $B_k=0$ results in a purely first-order
  algorithm.
  
  The condition \req{Bbound-a} on the Hessian approximations is quite weak,
  and allows in particular for a variety of quasi-Newton approaches,
  limited-memory or otherwise. In a finite-sum context, sampling bounded
  Hessians is also possible.
\item Conditions \req{GCPcond-a}--\req{gammak-def-a} define a ``generalized
    Cauchy point'' (GCP), much in the spirit of standard trust-region methodology
    (see \cite[Section~6.3]{ConnGoulToin00} for instance), where the
    quadratic model \req{q-model} is minimized in \req{gammak-def-a} along a
    good first-order direction ($s_k^L$) to obtain a ``Cauchy step''
    $s_k^Q$. Any step $s_k$ can then be accepted  provided it remains in the
    trust region (see \req{sbound-a}) and enforces a decrease in the quadratic
    model which is a least a fraction $\tau$ of that achieved at the Cauchy
    step (see \req{GCPcond-a}). 
\item At variance with many existing trust-region algorithms, the radius
  $\Delta_k$ of the trust-region \req{Delta-def-a} is not recurred adaptively
  from iteration to iteration depending on how well the quadratic model
  predicts function values, but is directly defined as a scaled version of
  the local gradient. This is not without similarities with the trust-region
  methods proposed by \cite{FanYuan01}, which corresponds to a scaling
  factor equal to $\|g_k\|^{-1}$, or \cite{CurtScheShi19} where 
  the trust-region radius depends on $\|g_k\|$.
\item As stated, the \al{ASTR1} algorithm does not include a termination
  rule, but such a rule can easily be introduced by terminating the algorithm in
  Step~1 if $\|g_k\| \leq \epsilon$, where $\epsilon>0$ is a user-defined
  first-order accuracy threshold.
\item It may seem to the reader that we have introduced two
  algorithmic parameters typically not present in existing OFFO
  methods. As it turns out, this is standard practice for
  trust-region methods and it is widely acknowledged that the
  behaviour of the algorithm is relatively insensitive to the choice
  made. Typically value are
  \[
  \tau = \sfrac{1}{10} \tim{ and } \kB= 10^5,
  \]
  the last one being possibly adapted to reflect the problem scaling.
  Note that these values are constant throughout the execution of the
  algorithm.  At variance, $\gamma_k$ is the iteration dependent
  stepsize, a quantity present in every first-order minimization
  method. Observe that we do not impose restrictions of the stepsize
  (beyond being positive), thereby covering most standard choices.
  Note that $\gamma_k=1$ and $s_k^Q=s_k^L$ whenever $B_k=0$.
\end{enumerate}

\noindent
The algorithm being defined, the first step of our analysis is to derive a
fundamental property of objective-function decrease, valid for all choices of
the scaling factors satisfying AS.4.

\llem{lemma:GCP}{
Suppose that AS.1, AS.2 and AS.4 hold. Then we have that, for all
$k\ge0$,
\beqn{gen-decr}
f(x_{j+1})
\le f(x_j) -\sum_{i=1}^n \frac{\tau \varsigma_{\min} g_{i,j}^2}{2\kB w_{i,j}}
      + \half(\kB + L) \sum_{i=1}^n \frac{g_{i,j}^2}{w_{i,j}^2}
\eeqn
and
\beqn{fdecrease}
f(x_0)-f(x_{k+1})
\ge \sum_{j=0}^k \sum_{i=1}^n\frac{g_{i,j}^2}{2\kB w_{i,j}}
    \left[\tau\varsigma_{\min}-\frac{\kBBL}{w_{i,j}}\right]
\eeqn
where $\varsigma_{\min} \eqdef \min_{i\in\ii{n}}\varsigma_i$ and $\kBBL \eqdef \kB(\kB+L)$.
}

\proof{
Using \req{sL-def-a} and AS.4, we deduce that, for every $j\geq 0$,
\beqn{gs1-a}
|g_j^T s_j^L|
= \sum_{i=1}^n \frac{g_{i,j}^2}{w_{i,j}}
= \sum_{i=1}^n \frac{w_{i,j}g_{i,j}^2}{w_{i,j}^2}
\ge \sum_{i=1}^n \frac{\varsigma_ig_{i,j}^2}{w_{i,j}^2}
\ge \varsigma_{\min}\|s_j^L\|^2.
\eeqn
Suppose first that $(s_j^L)^TB_js_j^L>0$ and $\gamma_j < 1 $.  Then, in view of \req{sQ-def-a},
\req{gammak-def-a}, \req{gs1-a} and \req{Bbound-a},
\[
g_j^Ts_j^Q+\half (s_j^Q)^TB_js_j^Q
= \gamma_j g_j^Ts_j^L+\half \gamma_j^2 (s_j^L)^TB_js_j^L
= -\frac{(g_j^Ts_j^L)^2}{2(s_j^L)^TB_js_j^L}
\le - \frac{\varsigma_{\min}|g_j^Ts_j^L|}{2\kB}.
\]
Combining this inequality with the first equality in \req{gs1-a} then gives that
\beqn{mdecrease-a}
g_j^Ts_j^Q+\half (s_j^Q)^TB_js_j^Q
\le - \frac{\varsigma_{\min}}{2\kB}\sum_{i=1}^n \frac{g_{i,j}^2}{w_{i,j}}.
\eeqn
Suppose now that
$(s_j^L)^TB_js_j^L \le 0$ or $\gamma_j=1$. Then, using \req{sQ-def-a}, \req{mdecrease-a} and \req{sL-def-a},
\[
g_j^Ts_j^Q+\half (s_j^Q)^TB_js_j^Q
= g_j^Ts_j^L+\half (s_j^L)^TB_js_j^L
\le \half g_j^Ts_j^L < 0
\]
and \req{mdecrease-a} then again follows from the bound $\kB \ge
1$. Successively using AS.1--AS.2, \req{GCPcond-a}, \req{mdecrease-a},
\req{Bbound-a} and \req{Delta-def-a} then gives that, for $j\ge0$,
\begin{align*}
f(x_{j+1})
& \le f(x_j) + g_j^Ts_j + \half s_j^TB_js_j - \half s_j^TB_js_j + \half L\|s_j\|^2 \label{posdef-a}\\*[1.5ex]
& \le f(x_j) + \tau\left(g_j^Ts_j^Q + \half (s_j^Q)^TB_js_j^Q\right)  + \half(\kB + L) \|s_j\|^2 \\
& \le f(x_j) -\sum_{i=1}^n \frac{\tau \varsigma_{\min} g_{i,j}^2}{2\kB w_{i,j}}
      + \half(\kB + L) \sum_{i=1}^n \Delta_{i,j}^2 \\
& \le f(x_j) -\sum_{i=1}^n \frac{\tau \varsigma_{\min} g_{i,j}^2}{2\kB w_{i,j}}
      + \half(\kB + L) \sum_{i=1}^n \frac{g_{i,j}^2}{w_{i,j}^2} 
\end{align*}
This is  \req{gen-decr}.
Summing up this inequality for $j\in\iiz{k}$ then yields \req{fdecrease}.
} 

\noindent
Armed with Lemma~\ref{lemma:GCP}, we are now in position to specify particular
choices of the scaling factors $w_{i,k}$ and derive the convergence properties
of the resulting variants of \al{ASTR1}.

\numsection{An Adagrad-like variant of \tal{ASTR1} using second-order models}\label{adap-s}

We first consider a choice of scaling factors directly derived from the
definition of the Adagrad algorithm \cite{DuchHazaSing11}.
For given
$\varsigma \in (0,1]$, $\vartheta \in (0,1]$, $\theta > 0$ and $\mu \in (0,1)$, define,
for all $i\in\ii{n}$ and for all $k \geq 0$,
\beqn{w-adag}
w_{i,k} \in \left[\sqrt{\vartheta}\, v_{i,k}, v_{i,k}\right]
\tim{ where }
v_{i,k} \eqdef \theta\left(\varsigma + \sum_{\ell=0}^k g_{i,\ell}^2\right)^\mu.
\eeqn
The Adagrad scaling factors are recovered by $\mu = \half$, $\theta = 1$ and $\vartheta=1$, and
\al{ASTR1} with \req{w-adag} and $B_k=0$ is then the standard (deterministic) Adagrad
method. The formulation \req{w-adag}  allows a parametric analysis of methods ``in
the neighbourhood'' of Adagrad, using not only first-order but also second-order information.
The $\vartheta$ parameter is introduced for flexibility, in
particular allowing non-monotone scaling factors\footnote{Typical
values are $\varsigma_i = \sfrac{1}{100}$ and $\vartheta = \sfrac{1}{1000}$.}
The additional scaling parameter\footnote{Which can be viewed as 
a stepsize/learning rate parameter when $B_k=0$.} $\theta$ is introduced 
as a technique to improve the convergence rate of the resulting algorithm.  
Such a scaling may be useful when the gradient is sparse
\cite{DuchJordMcMa13} or when designing a private Adagrad version
\cite{Asietal21} in order to improve the complexity bound with respect
to the problem's parameters. The above parametrization  
with $\theta=1$ has also been considered in \cite{ChakChop21} in the
context of a continuous Ordinary Differential Equations analysis and,
with $\theta$ depending on problem's constants, in \cite{LiOrab19}
where the sum on $\ell$ in \req{w-adag} is terminated at $\ell = k-1$
and $\mu$ is restricted to the interval $[\frac{1}{2}, 1) $ in the
  stochastic regime. In what follows, we consider the discrete
  deterministic case for the interval $(0,1)$. 

Before stating the global rate of convergence of the variant of \al{ASTR1}
using \req{w-adag}, we first prove a lemma, partly inspired by
\cite{WardWuBott19,DefoBottBachUsun22}.

\llem{gen:series}{Let $\{a_k\}_{k\ge 0}$ be a non-negative sequence,
$\alpha > 0$, $\xi>0$ and define, for each $k \geq 0$,
$b_k = \sum_{j=0}^k a_j$.  Then if $\alpha \neq 1 $,
\beqn{allalpha series-bound}
\sum_{j=0}^k  \frac{a_j}{(\xi+b_j)^{\alpha}}
\le \frac{1}{(1-\alpha)} ( (\xi + b_k)^{1-  \alpha} - \xi^{1-  \alpha} ).
\eeqn
Otherwise (i.e.\ if $\alpha  = 1$),
\beqn{alphasup1series-bound}
\sum_{j=0}^k  \frac{a_j}{(\xi+b_j)}
\le  \log\left(\frac{\xi + b_k}{\xi} \right).
\eeqn
}

\proof{
Consider first the case where  $\alpha \neq 1$ and note that
$\frac{1}{(1-\alpha)} x^{1-\alpha}$ is then a non-decreasing and concave
function on $(0,+\infty)$. Setting $b_{-1} = 0$ and using these properties, we
obtain that, for $j\geq 0$, 
\begin{align*}
\frac{a_j}{(\xi + b_j)^\alpha} &\leq \frac{1}{1-\alpha}\left( (\xi
  + b_j)^{1-\alpha}- (\xi + b_j - a_j)^{1-\alpha}\right) \\
&= \frac{1}{1-\alpha} \left( (\xi + b_j)^{1-\alpha} - (\xi + b_{j-1})^{1-\alpha}\right).
\end{align*} 
We then obtain \eqref{allalpha series-bound} by summing this inequality for $j\in\iiz{k}$.

Suppose now that $\alpha = 1$, we then use the concavity and non-decreasing
character of the logarithm to derive that
\[
\frac{a_j}{(\xi + b_j)^\alpha}
= \frac{a_j}{(\xi + b_j)}
\leq  \log(\xi + b_j) - \log(\xi + b_j - a_j) 
=  \log(\xi + b_j) - \log(\xi + b_{j-1}).
\]
The inequality \req{alphasup1series-bound} then again follows by summing for $j\in\iiz{k}$.
}

\noindent
From \eqref{allalpha series-bound}, we also obtain that, for $\alpha < 1$,
\beqn{alphainf1maj}
\sum_{j=0}^k \frac{a_j}{(\xi+b_j)^{\alpha}} \le \frac{1}{(1-\alpha)}(\xi + b_k)^{1-\alpha}
\eeqn
while, for $\alpha > 1$, 
\beqn{alphasup1maj}
\sum_{j=0}^k \frac{a_j}{(\xi+b_j)^{\alpha}} \le \frac{\xi^{1-\alpha}}{(\alpha-1)}.
\eeqn

\noindent
Note that both the numerator and the denominator of the right-hand side of
\req{allalpha series-bound} tend to zero when $\alpha$ tends to one. Applying
l'Hospital rule, we then see that this right-hand side tends to the right-hand
side of \req{alphasup1series-bound} and the bounds on $\sum_{j=0}^k
a_j/(\xi+b_j)^{\alpha}$ are therefore continuous at $\alpha = 1$.

Lemma~\ref{gen:series} is crucial in the proof of our main complexity
result, which we now state.

\lthm{theorem:allmu}{Suppose that AS.1--AS.3 hold and that the
\al{ASTR1} algorithm is applied to problem \req{problem} with its
scaling given by \req{w-adag}.
If we define
\[
\Gamma_0 \eqdef f(x_0)-\flow,
\]
then,
\begin{itemize}
\item[(i) ] if $ 0< \mu < \half$,
\beqn{gradboundmuinfhalf}
\average_{j\in\iiz{k}}\|g_j\|^2 \le  \frac{\kappa_1}{k+1},
\eeqn
with
\beqn{k3-def}
\kappa_1 = \max\left\{
           \varsigma,
           \left[\frac{2^{2\mu}\vartheta(1-2\mu) \theta^2\Gamma_0}{n(\kB+L)}\right]^{\sfrac{1}{1-2\mu}},\left[\frac{4\,n\,\kBBL}{(1-2\mu)\tau \theta \varsigma^\mu\vartheta^\sfrac{3}{2}}\right]^\sfrac{1}{\mu} 
           \right\};
\eeqn
\item[(ii) ] if $ \mu = \half$, 
\beqn{gradbound}
\average_{j\in\iiz{k}}\|g_j\|^2
\le \frac{\kappa_2}{k+1},
\eeqn
with \beqn{k6-def}
\kappa_2 = \max\left\{
\varsigma,\bigfrac{1}{2}e^{\frac{2 \Gamma_0\vartheta \theta^2}{n(\kappa_B + L)}},
\frac{1}{2\varsigma}
\left(\frac{8n\kB(\kB+L)}{\tau\vartheta^\sfrac{3}{2} \theta}\right)^2
\,\left|W_{-1}\left(-\frac{\tau \varsigma\theta\vartheta^\sfrac{3}{2}}{8n\kB(\kB+L)}\right)\right|^2
\right\},
\eeqn
where $W_{-1}$ is the second branch of the Lambert function \cite{Corletal96};
\item[(iii) ] if $ \half< \mu < 1$,
\beqn{gradboundmusupphalf}
\average_{j\in\iiz{k}}\|g_j\|^2
\le \frac{\kappa_3}{k+1}
\eeqn
with 
\beqn{k7-def}
\hspace*{-4mm}\kappa_3 \eqdef
= \max\left\{\varsigma,\left[ \frac{2^{1+\mu} \kB}{\tau \varsigma^\mu\sqrt{\vartheta}}
  \left(\!\Gamma_0 \theta+\frac{n(\kB+L)\varsigma^{1-2\mu}}{2\vartheta \theta(2\mu-1)}\right)\right]^{\sfrac{1}{1-\mu}}\right\}
\eeqn
\end{itemize}
}

\proof{
We see from \req{w-adag} that $w_{i,k}$ verifies \textbf{AS.4}. We may
thus use Lemma~\ref{lemma:GCP}. Moreover, \req{w-adag} also implies that
\beqn{eq3.1}
\varsigma^\mu \sqrt{\vartheta} \theta
\leq w_{i,j}
\leq \theta \left(\varsigma + \sum_{\ell=0}^j\|g_\ell\|^2\right)^\mu
\eeqn
for all $j\geq0$ and all $i\in\ii{n}$.
We now deduce from
\req{Delta-def-a} and \req{fdecrease} that, for $k\geq 0$, 
\beqn{summing-up1}
\begin{array}{lcl}
f(x_{k+1})
& \le & f(x_0) - \bigsum_{j=0}^k \bigfrac{\tau \varsigma^\mu\sqrt{\vartheta}\,\|g_j\|^2}{2 \kB \max_{i\in\ii{n}} w_{i,k}}
+  \half (\kB + L) \bigsum_{i=1}^n \bigsum_{j=0}^k \Delta_{i,j}^2.
\end{array}
\eeqn
For each $i\in\ii{n}$, we then apply Lemma~\ref{gen:series}
with $a_\ell = g_{i,\ell}^2$, $\xi=\varsigma$ and $\alpha= 2 \mu < 1$, and obtain
from \req{Delta-def-a} and \req{w-adag} that,
\beqn{harsh-bound}
\bigsum_{j=0}^k \Delta_{i,j}^2
\le \bigfrac{1}{\theta^2\vartheta(1-2 \mu)} \left[\left(\varsigma + \bigsum_{\ell=0}^k g_{i,\ell}^2\right)^{1-2\mu}
  - \varsigma^{1-2 \mu} \right]
\le  \bigfrac{1}{\theta^2\vartheta(1-2 \mu)} \left(\varsigma + \bigsum_{\ell=0}^k
g_{i,\ell}^2\right)^{1-2\mu}.
\eeqn
Now
\beqn{ineq2}
\bigsum_{i=1}^n \bigsum_{j=0}^k \Delta_{i,j}^2
\le  \bigsum_{i=1}^n \bigfrac{1}{\theta^2\vartheta(1-2 \mu)} \left(\varsigma + \bigsum_{i=1}^n \bigsum_{\ell=0}^k
g_{i,\ell}^2\right)^{1-2\mu}
\le  \bigfrac{n}{\theta^2\vartheta(1-2 \mu)}\left(\varsigma + \bigsum_{\ell=0}^k \|g_\ell\|^2\right)^{1-2\mu}
\eeqn
and thus substituting this bound in \req{summing-up1} and using AS.3 gives that
\beqn{eq01}
\sum_{j=0}^k \bigfrac{\tau \varsigma^\mu \sqrt{\vartheta}\|g_j\|^2}{2 \kB \max_{i\in\ii{n}} w_{i,k}}
\le \Gamma_0 + \frac{n(\kB + L)}{2 \theta^2 \vartheta}(1-2\mu)\left(\varsigma + \sum_{j=0}^k \|g_j\|^2\right)^{1-2\mu}.
\eeqn
Suppose now that
\beqn{hypbad}
\sum_{j=0}^k \|g_j\|^2
\geq \max\left\{ \varsigma,\left[\frac{2^{2\mu} \theta^2\vartheta(1-2\mu)\Gamma_0}{n(\kB+L)}\right]^{\sfrac{1}{1-2\mu}} \right\},
\eeqn
implying
\[
\varsigma + \sum_{j=0}^k \|g_j\|^2 \leq 2\sum_{j=0}^k \|g_j\|^2
\tim{and}
\Gamma_0 \leq \frac{n(\kB + L)}{2 \theta^2\vartheta(1-2\mu)}\left(2\sum_{j=0}^k \|g_j\|^2\right)^{1-2\mu}.
\]
Then, using \req{eq01} and \req{eq3.1},
\[
\frac{\tau \varsigma^\mu\sqrt{\vartheta}}{2^{1+\mu} \,\kB  \theta  \left[ \sum_{\ell=0}^k \|g_\ell\|^2\right]^\mu}\sum_{j=0}^k \|g_j\|^2
\leq \frac{2^{1-2\mu}\,n\,(\kB+L)}{\vartheta \theta^2(1-2\mu)} \left(\sum_{j=0}^k \|g_j\|^2\right)^{1-2\mu}.
\]
Solving this inequality for $\sum_{j=0}^k \|g_j\|^2$ gives that
\[
\sum_{j=0}^k \|g_j\|^2 \leq 
\left[\frac{ 4\,n\,\,\kBBL}{(1-2\mu) \theta\tau \varsigma^\mu\vartheta^\sfrac{3}{2}}\right]^{\sfrac{1}{\mu}}
\]
and therefore
\beqn{ab1}
\average_{j\in \iiz{k}} \|g_j\|^2
\leq \left[\frac{4\,n\,\kBBL}{(1-2\mu) \theta\tau \varsigma^\mu\vartheta^\sfrac{3}{2}}\right]^{\sfrac{1}{\mu}}\cdot\frac{1}{{k+1}}.
\eeqn
Alternatively, if \req{hypbad} fails, then
\beqn{ab2}
\average_{j\in\iiz{k}} \|g_j\|^2
< \max\left\{
\varsigma,\left[\frac{2^{2\mu}\vartheta(1-2\mu) \theta^2\Gamma_0}{2^{1-2\mu}n(\kB+L)}\right]^{\sfrac{1}{1-2\mu}}
\right\}\cdot\frac{1}{k+1}.
\eeqn
Combining \req{ab1} and \req{ab2}
gives \req{gradboundmuinfhalf}.
	
Let us now consider the case where $\mu = \half$.
For each $i\in\ii{n}$, we apply Lemma~\ref{gen:series}
with $a_k = g_{i,k}^2$, $\xi = \varsigma$ and $\alpha= 2 \mu = 1$ and obtain that,
\[
\bigsum_{i=1}^n\bigsum_{j=0}^k \Delta_{i,j}^2
\le \frac{1}{\vartheta \theta^2}\bigsum_{i=1}^n
\log\left(\frac{1}{\varsigma}\left(\varsigma+ \bigsum_{\ell=0}^k
g_{i,\ell}^2\right)\right)
\le \frac{n}{\vartheta \theta^2} \log\left(1 +  \frac{1}{\varsigma}\bigsum_{\ell=0}^k \| g_\ell\|^2\right).
\]
and substituting this bound in \eqref{summing-up1} then gives that
\[
\sum_{j=0}^k \bigfrac{\tau \sqrt{\varsigma\vartheta} \|g_j\|^2}{2 \kB \max_{i\in\ii{n}} w_{i,k}}
\le \Gamma_0 + \frac{n(\kB + L)}{2\vartheta \theta^2} \log\left( 1 +  \frac{1}{\varsigma}\bigsum_{j=0}^k \| g_j\|^2\right).
\]
Suppose now that 
\beqn{logmuhalfhyp}
\sum_{j=0}^k \|g_j\|^2
\geq \max\left[ \varsigma, \bigfrac{1}{2}e^{\frac{2 \vartheta \theta^2 \Gamma_0}{n(\kappa_B + L)}}\right],
\eeqn
\noindent
implying that 
\[
1 + \frac{1}{\varsigma}\sum_{j=0}^k \|g_j\|^2 \leq \frac{2}{\varsigma}\sum_{j=0}^k \|g_j\|^2
\tim{and}
\Gamma_0 \leq \frac{n(\kB + L)}{2\vartheta \theta^2} \log\left( \frac{2}{\varsigma}\bigsum_{j=0}^k \| g_j\|^2\right).
\]
\noindent
Using \req{eq3.1} for $\mu=\half$, we obtain then that
\[
\bigfrac{\tau \sqrt{\varsigma\vartheta} }{2\sqrt{2}\, \theta\, \kB \,\sqrt{\sum_{\ell=0}^k \| g_\ell\|^2}} \bigsum_{j=0}^k  \|g_j\|^2
\le \frac{n(\kB + L)}{\vartheta \theta^2} \log\left(  \frac{2}{\varsigma} \bigsum_{j=0}^k \|
g_j\|^2 \right),
\]
that is
\beqn{tosolvmuhalf}
 \bigfrac{\tau \sqrt{2\varsigma} \vartheta^\sfrac{3}{2} \theta}{4 \kB } \sqrt{  \bigsum_{j=0}^k \| g_j\|^2}
\le 2 n(\kB + L) \log\left(  \sqrt{\frac{2}{\varsigma} \bigsum_{j=0}^k \| g_j\|^2}
\right).
\eeqn
Now define
\beqn{ggu}
\gamma_1 \eqdef \bigfrac{\tau \varsigma\vartheta^\sfrac{3}{2} \theta }{4 \kB },
\ms
\gamma_2 \eqdef  2 n(\kB + L)
\tim{ and }
u \eqdef \sqrt{\frac{2}{\varsigma} \bigsum_{j=0}^k\|g_j\|^2}
\eeqn
and observe that that $\gamma_2 > 3 \gamma_1$ because $\tau\sqrt{\varsigma} \vartheta^\sfrac{3}{2}\leq 1$
and $\kB \geq 1$.
The inequality \req{tosolvmuhalf} can then be rewritten as
\beqn{tosolvemuhalf}
\gamma_1 u \le \gamma_2 \log(u).
\eeqn
Let us denote by $\psi(u) \eqdef \gamma_1 u - \gamma_2 \log(u)$. Since
$\gamma_2 > 3 \gamma_1$, the equation $\psi(u)=0$ admits two roots $u_1 \leq
u_2$ and \req{tosolvemuhalf} holds for $u\in[u_1,u_2]$.
The definition of $u_2$ then gives that
\[
\log(u_2)- \frac{\gamma_1}{\gamma_2}u_2 = 0
\]
which is
\[
u_2e^{-\frac{\gamma_1}{\gamma_2}u_2} = 1.
\]
Setting $z = -\frac{\gamma_1}{\gamma_2}u_2$, we obtain that
\[
z e^z = -\frac{\gamma_1}{\gamma_2}
\]
Thus $z = W_{-1}(-\frac{\gamma_1}{\gamma_2})<0$, where $W_{-1}$ is the second
branch of the Lambert function defined over $[-\frac{1}{e}, 0)$.
As $-\frac{\gamma_1}{\gamma_2} \geq -\frac{1}{3} $, $z$ is well defined and thus
\[
u_2
= -\frac{\gamma_2}{\gamma_1}\,z
= -\frac{\gamma_2}{\gamma_1}\,W_{-1}\left(-\frac{\gamma_1}{\gamma_2}\right)>0.
\]
As a consequence, we deduce from \req{tosolvemuhalf} and \req{ggu} that
\[
\bigsum_{j=0}^k\|g_j\|^2
= \frac{\varsigma}{2}\,u_2^2
=\frac{1}{2\varsigma}\,\left(\frac{8n\kB(\kB+L)}{\tau\vartheta^\sfrac{3}{2} \theta}\right)^2
\,\left|W_{-1}\left(-\frac{\tau \varsigma\vartheta^\sfrac{3}{2} \theta}{8n \kB(\kB+L)}\right)\right|^2.
\]
and
\beqn{firstmuhalfineq}
\average_{j\in\iiz{k}} \|g_j\|^2
\leq \frac{1}{2\varsigma}\,
\left(\frac{8n\kB(\kB+L)}{\tau\vartheta^\sfrac{3}{2} \theta}\right)^2
\,\left|W_{-1}\left(-\frac{\tau \varsigma \vartheta^\sfrac{3}{2} \theta}{8n
  \kB(\kB+L)}\right)\right|^2\cdot\frac{1}{k+1}.
\eeqn
If \eqref{logmuhalfhyp} does not hold, we have that
\beqn{secmuhalfineq}
\average_{j\in\iiz{k}} \|g_j\|^2
< \max\left\{
\varsigma, \frac{1}{2}\, e^{\frac{2\Gamma_0\vartheta \theta^2}{n(\kappa_B + L)}}\right\}\cdot \frac{1}{k+1}.
\eeqn
Combining \req{firstmuhalfineq} and \req{secmuhalfineq}
gives \req{gradbound}.

Finally, suppose that $\half < \mu < 1$.
Once more, we apply Lemma~\ref{gen:series} for each $i\in\ii{n}$
with $a_\ell = g_{i,\ell}^2$, $\xi=\varsigma$ and $\alpha= 2 \mu > 1$ and obtain that
\beqn{mugeqhalfound}
\sum_{j=1}^k \Delta_{k,j}^2
\le \frac{1}{\theta^2\vartheta(1-2 \mu)} \left( \Big(\varsigma + \sum_{\ell=0}^k g_{i,\ell}^2\Big)^{1-2\mu}- \varsigma^{1-2\mu} \right)
\le \frac{\varsigma^{1-2 \mu}}{\vartheta \theta^2 (2\mu -1)}.
\eeqn
Substituting the bound \req{mugeqhalfound} in \req{summing-up1} and using
\req{eq3.1} and AS.3 gives that
\[
\sum_{j=0}^k  \frac{1}{(\varsigma +
	\sum_{j=0}^k \| g_j\|^2 )^\mu}\bigfrac{\tau \varsigma^\mu\sqrt{\vartheta}\|g_j\|^2}{2 \kB \theta}
\le \Gamma_0 + \frac{n(\kB + L)\varsigma^{1-2\mu}}{2 \theta^2\vartheta(2\mu-1)}. \\
\]
If we now suppose that
\beqn{musuphalfhyp}
\sum_{j=0}^k \| g_j\|^2 \geq \varsigma,
\eeqn
then
\beqn{firstmugeqhalf}
\average_{j\in\iiz{k}} \|g_j\|^2
\le \left[ \frac{2^{1+\mu} \kB}{\tau \varsigma^\mu\sqrt{\vartheta}} \left(\Gamma_0 \theta
  + \frac{n(\kB + L)\varsigma^{1-2\mu}}{2\vartheta(2\mu-1) \theta}\right) \right]^{\sfrac{1}{1-\mu}}
  \cdot \frac{1}{k+1}.
\eeqn
If \req{musuphalfhyp} does not hold, we derive that
\beqn{secmugeqhalf}
\average_{j\in\iiz{k}} \|g_j\|^2 \leq \frac{\varsigma}{(k+1)}.
\eeqn
Thus, \req{firstmugeqhalf} and\req{secmugeqhalf}
finally imply \req{gradboundmusupphalf}. 
} 

\noindent
These result suggest additional remarks.
\begin{enumerate}

\item That the bounds given  are not continuous as a function $\mu$
   at $\mu = \half$ is a result of our bounding process within the
  proof of Theorem~\ref{theorem:allmu} (for instance in the last inequality of
  \req{harsh-bound}).  Continuous bounds have been proved
  (see\cite{GratJeraToin22b}) if one is ready to
  assume that the objective functions' gradients remain uniformly
  bounded. 
  
\item If the algorithm is terminated as soon as $\|g_k\| \le \epsilon$ (which
  is customary for deterministic algorithms searching for first-order points),
  it must stop at the latest at iteration
  \beqn{eps-order}
  k = \kappa_\star^2  \epsilon^{-2},
  \eeqn
  where $\kappa_\star=\kappa_1$ for $\mu\in(0,\half)$,
  $\kappa_\star=\kappa_2$ for $\mu=\half$ and $\kappa_\star=\kappa_3$ for
  $\mu\in(\half,1)$. It is truly remarkable that there exist first-order
  OFFO methods whose global complexity order is identical to that of
  standard first-order methods using function  evaluations (see
  \cite{Nest04,GratSartToin08,CartGoulToin10a} or
  \cite[Chapter~2]{CartGoulToin22}), despite the fact that the latter
  exploit  significantly more information. As mentionned in the
  introduction,  this complexity rate was also derived for the
  analysis of Adagrad-Norm in \cite{WardWuBott19}. 

\item
 It is possible to give a weaker but more explicit bound on $\kappa_2$ by finding an upper
 bound on the value of the involved Lambert function. This can be obtained by 
 using \cite[Theorem~1]{Chat13} which states that, for $x>0$,
 \beqn{Lamb-bound}
 \left|W_{-1}(-e^{-x-1})\right| \leq 1+ \sqrt{2x} + x.
 \eeqn
 Remembering that, for $\gamma_1$ and $\gamma_2$ given by \req{ggu},
 $\log\left( \frac{\gamma_2}{\gamma_1}\right) \geq \log(3) > 1$
 and taking $x = \log\left( \frac{\gamma_2}{\gamma_1}\right) -1 > 0$
 in \req{Lamb-bound} then gives that
 \[
 \left|W_{-1}\left(-\frac{\gamma_1}{\gamma_2} \right)\right|
 \le \log\left(\frac{\gamma_2}
     {\gamma_1}\right)+\sqrt{2\left(\log\left(\frac{\gamma_2}{\gamma_1}\right)-1\right)}.
 \]
 
\item It is also possible to extend the definition of $s_k^L$ in \req{sL-def-a}
  by premultiplying it by a stepsize $\alpha_k\in[\alpha_{\min}, 1]$ for some
  $\alpha_{\min} \in (0,1]$. Our results again remain valid (with modified
  constants). Covering a deterministic momentum-less Adam would require extending the
  results to allow for \req{w-adag} to be replaced by
  \beqn{w-adam}
  w_{i,k} = \varsigma + \sum_{j=0}^k\beta_2^{k-j}g_{i,j}^2  \ms (i \in \ii{n})
  \eeqn
  for some $\beta_2 < 1$.  This can be done by following the argument of Theorem~2
  in \cite{DefoBottBachUsun22}. However, as in this reference, the final bound
  on the squared gradient norms does not tend to zero when $k$ grows\footnote{A
  constant term in $-\log(\beta_2)$ refuses to vanish.}, illustrating the
  (known) lack of convergence of Adam. We therefore do not investigate this
  option in detail.
\item
   Focusing on the choice of the parameter $\mu$ independently of $\vartheta$ and 
   $\theta$, we verify that the choice $\mu = \frac{1}{2}$  is best, yielding an 
   upper complexity bound for the deterministic Adagrad algorithm and recovering 
   a similar result obtained in \cite{WardWuBott19} for the Adagrad-Norm algorithm.

   The choice $\mu = 1$ is not covered by our theory but has been considered in 
   \cite{MukkHein17},  where a specific choice of the constant $\theta$ with $\mu = 1$ 
   is used to derive a first order method with optimal regret for strongly convex 
   problems. Unfortunately, the proposed SC(Strong Convex)-Adagrad algorithm requires 
   the knowledge of the problem's Lipschitz constant and noise characteristics.
\end{enumerate}

We now discuss the dependence of the bounds given by Theorem~\ref{theorem:allmu} as 
a function of the problem's constants $L$, $n$ and $\Gamma_0$ and recall that 
if $n$ is known \textit{a priori}, this is generally not the case for the Lispchitz 
constant $L$, while the gap $\Gamma_0$ may be available for some classes of problems 
(such as nonlinear regressions where $\Gamma_0\leq f(x_0)$). We are mostly interested 
in the explicit dependence of the bounds on $n$, taking into account that the unknown 
$L$ often varies little with dimension --such as in problems arising from 
discretizations-- but admittedly ignoring the fact it may also depend on $n$, sometimes 
severely \cite[p.~14]{CartGoulToin22}.  As we are allowed to do so, we would like 
to choose the scaling parameter $\theta$ in order to offset this dependence as much as possible.
In view of \eqref{k3-def}, \eqref{k6-def} and \eqref{k7-def}, and noting that 
\req{Lamb-bound} indicates that the $|W_{-1}|^2$ term  in \eqref{k6-def} can be 
summarized as $\calO\left( \log(n L)^2\right)$, we attempt to balance the impact 
of the various terms involving both $\theta$ and $n$ and, in the absence of additional 
information on the problem, select
\beqn{thetabest}
\theta_\star = \left\{\begin{array}{ll}
\sqrt{n/\Gamma_0} & \tim{when $\Gamma_0$ is known,} \\
\sqrt{n} & \tim{otherwise.}
\end{array}\right.
\eeqn 
We may then compare the bounds obtained by reinjecting this value in  \eqref{k6-def} 
($\mu = \half$) with the results for deterministic Adagrad-Norm \cite{WardWuBott19} and 
Adagrad \cite{TraoPauw21}.  This comparison is  summarized in the following table.
\begin{center}
{\renewcommand{\arraystretch}{1.8}%
\begin{tabular}{|l|c|c|}
\hline
Regime & Algorithm & Constant dependence  \\
\hline
Non-Convex \cite{WardWuBott19} & $w_{i,k} = \sqrt{\varsigma + \sum_{j= 0}^k \|g_j \|^2}$   &  $\mathcal{O}\left( L^2 \log(L)^2 \right)$  \\
\hline
Convex  \cite{TraoPauw21} &  \eqref{w-adag} with ($\theta =1, \, \vartheta = 1$)   & $\mathcal{O}\left( n^2L^2 \log(L)^2 \right)$\\
\hline
Non-convex \eqref{k6-def} & \eqref{w-adag} with ($\theta =\sqrt{n}, \, \vartheta = 1$) & $\mathcal{O}\left(n L^2 \log(\sqrt{n}L)^2  \right)$\\
\hline
\end{tabular}}
\end{center}

\noindent
Are the (good) upper bound given Theorem~\ref{theorem:allmu} sharp?

\lthm{sharp1}{The bounds \req{gradboundmuinfhalf}, \req{gradbound} and
  \req{gradboundmusupphalf} are essentially sharp in that, for each
  $\mu\in(0,1)$ and each $\eta\in(0,1]$, there exists a univariate function $f_{\mu,\eta}$
  satisfying AS.1-AS.3 such that, when applied to minimize $f_{\mu,\eta}$ from the
  origin,  the \al{ASTR1} algorithm  with \req{w-adag}, $B_k=0$ and $\vartheta=\theta=1$ produces a sequence of
  gradient norms given by $\|g_k\| = \frac{1}{k^{\half+\eta}}$.
}

\proof{Following ideas of \cite[Theorem~2.2.3]{CartGoulToin22}, we first
construct a sequence of iterates  $\{x_k\}$ for which
$f_{\mu,\eta}(x_k) = f_k$ and $\nabla_x^1f_{\mu,\eta}(x_k) = g_k$
for associated sequences of function and gradient values
$\{f_k\}$ and $\{g_k\}$, and then apply Hermite interpolation
to exhibit the function $f_{\mu,\eta}$ itself.  We start by defining
\beqn{gk-def}
g_0 \eqdef -2,
\ms
g_k \eqdef -\frac{1}{k^{\half+\eta}} \ms(k>0),
\eeqn
\beqn{sk-def}
s_0 \eqdef \frac{2}{(\varsigma+4)^\mu},
\ms
s_k \eqdef \frac{1}{k^{\half+\eta}(\varsigma+\sum_{j=0}^k g_j^2)^\mu} \ms(k>0)
\eeqn
yielding that
\beqn{gksk}
|g_0s_0| = \frac{4}{(\varsigma+4)^\mu},
\ms
|g_ks_k|
= \frac{1}{k^{1+2\eta}(\varsigma+\sum_{j=0}^k g_j^2)^\mu}
\leq \frac{1}{k^{1+2\eta}} \ms(k>0)
\eeqn
(remember that $g_0^2=4$).
We then define $B_k \eqdef 0$ for all $k\geq0$,
\beqn{xk-def}
x_0 = 0,
\ms
x_{k+1} = x_k+s_k \ms(k>0)
\eeqn
and
\beqn{fk-def}
f_0 = \frac{4}{(\varsigma+4)^\mu} + \zeta(1+2\eta)
\tim{ and }
f_{k+1} = f_k + g_ks_k \ms (k \geq 0),
\eeqn
where $\zeta(\cdot)$ is the Riemann zeta function.
Observe that the sequence $\{f_k\}$ is decreasing and that, for all $k\geq 0$,
\[
f_{k+1}
= f_0 - \bigsum_{k=0}^k|g_ks_k|
\geq f_0 - \frac{4}{(\varsigma+4)^\mu}- \bigsum_{k=1}^k\frac{1}{k^{1+2\eta}}
\geq f_0 - \frac{4}{(\varsigma+4)^\mu} - \zeta(1+2\eta)
\]
where we used \req{fk-def} and \req{gksk}. Hence
\req{fk-def} implies that
\beqn{fk-bound}
f_k \in [0, f_0] \tim{for all} k\geq 0.
\eeqn
Also note that, using \req{fk-def},
\beqn{dfok}
|f_{k+1} - f_k - g_ks_k| = 0,
\eeqn
while, using \req{sk-def},
\[
|g_0-g_1| = 1 \leq \frac{1}{2}(\varsigma+4)^\mu\,s_0.
\]
Moreover, using the fact that $1/x^{\half+\eta}$ is a convex function of $x$ over
$[1,+\infty)$, and that from \req{sk-def}
$s_k \geq \frac{1}{k^{\half + \eta} \left(\varsigma + 4 + k\right)^\mu}$,  we
derive that, for $k>0$,
\begin{align*}
|g_{k+1} - g_k| &= \left| \frac{1}{(k+1)^{\half+\eta}} - \frac{1}{k^{\half+\eta}} \right| \\
&\leq \left(\frac{1}{2}+\eta\right) \frac{1}{k^{\sfrac{3}{2} + \eta}} \\
&\leq \frac{3}{2} \, \frac{(\varsigma + 4 + k)^\mu}{ k  k^{\half + \eta} (\varsigma + 4 + k)^\mu} \\
&\leq \frac{3}{2} \, \frac{(\varsigma + 4 + k)^\mu}{k } s_k \\
&\leq \frac{3}{2} \, \left(\varsigma + 5 \right)^\mu s_k .
\end{align*}
These last bounds and \req{fk-bound} allow us to use standard Hermite
interpolation on the data given by $\{f_k\}$ and $\{g_k\}$: see, for instance,
Theorem~A.9.1 in \cite{CartGoulToin22} with $p=1$ and
\beqn{unc-kappaf-slow}
\kappa_f = \max\left[\frac{3}{2}(\varsigma+5)^\mu, f_0,2\right]
\eeqn
(the second term in the max bounding $|f_k|$ because of
\req{fk-bound} and the third bounding $|g_k|$ because of \req{gk-def}).
We then deduce that there exists a continuously differentiable function $f_{\mu,\eta}$
from $\Re$ to $\Re$ with Lipschitz continuous gradient  (i.e. satisfying
AS.1 and AS.2) such that, for $k\geq 0$,
\[
f_{\mu,\eta}(x_k) = f_k\tim{ and } \nabla_x^1f_{\mu,\eta}(x_k) = g_k.
\]
Moreover, the range of $f_{\mu,\eta}$ and $\nabla_x^1 f_{\mu,\eta}$ are constant independent
of $\eta$, hence guaranteeing AS.3 and AS.3.
The definitions \req{gk-def}, \req{sk-def}, \req{xk-def} and \req{fk-def}
imply that the sequences $\{x_k\}$,  $\{f_k\}$ and $\{g_k\}$ can be seen as generated by the
\al{ASTR1} algorithm (with \req{w-adag}, $B_k=0$ and
$\vartheta=\theta=1$) applied to $f_{\mu,\eta}$, starting from $x_0=0$
and the desired conclusion follows.
} 

\noindent
The bounds \req{gradboundmuinfhalf}, \req{gradbound} and
  \req{gradboundmusupphalf} are therefore \emph{essentially sharp} (in the sense of
\cite{CartGoulToin18a}) for the \al{ASTR1} algorithm with \req{w-adag} and $\vartheta=\theta=1$, which
is to say that the lower complexity bound for the algorithm is arbitrarily close
to its upper bound.
Interestingly, the argument in the
proof of the above theorem fails for $\eta=0$, as this choice yields that
\[
\sum_{j=0}^k g_j^Ts_j \geq \sum_{j=0}^k \frac{1}{k(\varsigma+\log(k+1))^\mu}.
\]
Since
\[
\int_1^k \frac{dt}{t(\log(t+1))^\mu}
> \int_1^k \frac{dt}{(t+1)(\log(t+1))^\mu}
=\frac{(\log(k+1))^{1-\mu}}{1-\mu}-\frac{\log(2)^{1-\mu}}{1-\mu}
\]
tends to infinity as $k$ grows, this indicates (in view \req{fk-def}) that
AS.3 cannot hold.  Also note that \req{gk-def} implies that the
gradients remain uniformly bounded.


\numsection{A further ``diminishing stepsizes'' variation on this theme}\label{newclass-s}

We now use a different proof technique to design new variants of \al{ASTR1} with
a fast global $k$-order.  This is achieved by modifying the definition of the
scaling factors $w_{i,k}$, requiring them to satisfy a fairly general growth
condition explicitly depending on $k$, the iteration index.
More specifically, we will assume, in this section, that
the scaling factors $w_{i,k}$ are chosen such that, for some power
parameter $0< \nu \leq \mu < 1$, all $i\in\ii{n}$ and some constants
$\varsigma_i\in(0,1]$ and $\theta > 0$,
\beqn{wik-bound-ming}
\theta \max[\varsigma_i,v_{i,k}]\,  (k+1)^\nu \le w_{i,k} \le \theta \max[\varsigma_i,v_{i,k}]\,(k+1)^\mu
\ms
(k \geq 0),
\eeqn
where, for each $i$, the $v_{i,k}$ satisfy the properties that
\beqn{vikprop}
v_{i,k+1} > v_{i,k}  \tim{ implies that } v_{i,k+1} \leq |g_{i,k+1}|
\eeqn
and
\beqn{viklow}
v_{i,k} \geq |g_{i,k}| / h(k)
\eeqn
for some positive function  $h(k)$ only depending on $k$.
The motivation for introducing these new variants is the remarkable numerical
performance \cite{GratJeraToin22a} of particular choices where
\[
v_{i,k} = \max_{j\in\iiz{k}}|g_{i,j}|
\tim{ and }
v_{i,k} = \frac{1}{k+1}\sum_{j\in\iiz{k}}|g_{i,j}|
\]
which both satisfy \req{vikprop} and \req{viklow} (with $h(k)=1$ for
the first and $h(k)=k+1$ for the second).
We further illustrate this in Section~\ref{numerics-s}.

We start by proving a useful technical result.

\llem{lemma:divsprop}{Consider and arbitrary $i\in \ii{n}$ and suppose
that there exists a $j_\varsigma$ such that
\beqn{termissmall}
\min\left[\frac{g_{i,j}^2}{\varsigma_i},\frac{g_{i,j}^2}{v_{i,j}}\right]
\leq \varsigma_i
\tim{ for } j \geq j_\varsigma.
\eeqn
Then
\beqn{termlbound}
\min\left[\frac{g_{i,j}^2}{\varsigma_i},\frac{g_{i,j}^2}{v_{i,j}}\right]
\geq \frac{g_{i,j}^2}{2\varsigma_i}
\tim{ for } j \geq j_\varsigma.
\eeqn
}

\proof{Suppose that there exists a $j>j_\varsigma$
such that $v_{i,j} >2\varsigma_i$.  Assume, without
loss of generality that $j$ is the smallest such index.
Then $v_{i,j} > v_{i,j-1}$ and \req{vikprop} implies that
$|g_{i,j}| \geq v_{i,j} \geq 2\varsigma_i$. As a consequence,
\[
\min\left[\frac{g_{i,j}^2}{\varsigma_i},\frac{g_{i,j}^2}{v_{i,j}}\right]
  \geq \min[4\varsigma_i,2\varsigma_i]
> \varsigma_i,
\]
which contradicts \req{termissmall}.  Thus no such $j$ can exists and $v_{i,j} \le
2 \varsigma_i$ for all $j > j_*$ and \req{termlbound} follows.
}

\noindent
We are now in position to state our complexity result for the
\al{ASTR1} algorithm using weight defined by \req{wik-bound-ming},
\req{vikprop} and \req{viklow}.

\lthm{theorem:ming-style}{
Suppose that AS.1, AS.2 and AS.3 hold and that the
\al{ASTR1} algorithm is applied to problem \req{problem}, where
the scaling factors $w_{i,k}$ are chosen in accordance with
\req{wik-bound-ming}, \req{vikprop} and \req{viklow}.
Then, for any $\eta \in (0,\tau \varsigma_{\min})$ and
\beqn{jstar-ming}
j_\eta \eqdef
\left(\frac{\kB(\kB+L)}{\theta \varsigma_{\min}(\tau
  \varsigma_{\min}-\eta)}\right)^{\sfrac{1}{\nu}},
\eeqn
there exist a constant $\kappa_\diamond$, a subsequence
$\{k_\ell\}\subseteq \{k\}_{j_\eta+1}^\infty$ and an index
$k_\varsigma$ (where $\kappa_\diamond$ and $k_\varsigma$
only depend on the problem and the algorithmic constants) such that, for all $k_\ell \ge k_\varsigma$,
\beqn{ngkbound-ming} 
\min_{j\in\iiz{k_\ell}}\|g_j\|^2
\le \kappa_\diamond \frac{(k_\ell+1)^\mu}{k_\ell-j_\eta}
\le \frac{2\kappa_\diamond (j_\eta+1)}{k_\ell^{1-\mu}}.
\eeqn
}

\proof{
From \req{fdecrease} and AS.3, using $w_{\min,j}\eqdef \min_{i\in\ii{n}}w_{i,k}$ ensures that
\beqn{newb-a-ming}
\Gamma_0
\ge f(x_0)-f(x_{k+1})
\ge \sum_{j=0}^k \sum_{i=1}^n \frac{g_{i,j}^2}{2\kB  w_{i,j}}
    \left[\tau \varsigma_{\min}-\frac{\kBBL}{w_{\min,j} }\right].
\eeqn
Consider now an arbitrary $\eta \in (0,\tau \varsigma_{\min})$ and suppose first
that, for some $j$, 
\beqn{brackle-ming}
\left[\tau \varsigma_{\min}-\frac{\kBBL}{w_{\min,j}}\right] \le \eta,
\eeqn
i.e., using \req{wik-bound-ming},
\[
\theta\varsigma_{\min}\,\,j^\nu \le w_{\min,j} \le \frac{\kBBL}{\tau \varsigma_{\min}-\eta}.
\]
But this is impossible for $j > j_\eta$ for $j_\eta$ given by \req{jstar-ming},
and hence \req{brackle-ming} fails for all $j > j_\eta$.
As a consequence, we have that, for $k> j_\eta$,
\begin{align}
f(x_{j_\eta+1}) - f(x_k)
&\ge \eta \sum_{j=j_\eta+1}^k \sum_{i=1}^n \frac{g_{i,j}^2}{2\kB  w_{i,j}}\neol
&\ge \frac{\eta}{2\kB} \sum_{j=j_\eta+1}^k \sum_{i=1}^n\frac{g_{i,j}^2}{\max[\varsigma_i,v_{i,j}]\, \theta\, (j+1)^\mu}\neol
&\ge \frac{\eta}{2\kB(k+1)^\mu \theta}
     \sum_{j=j_\eta+1}^k \sum_{i=1}^n\min\left[\frac{g_{i,j}^2}{\varsigma_i},\frac{g_{i,j}^2}{v_{i,j}}\right]\neol
&\ge \frac{\eta (k-j_\eta)}{2\kB(k+1)^\mu \theta}
     \min_{j\in\iibe{j_\eta+1}{k}}
     \left(\sum_{i=1}^n\min\left[\frac{g_{i,j}^2}{\varsigma_i},\frac{g_{i,j}^2}{v_{i,j}}\right]\right)\label{A1}
\end{align}
But we also know from \req{gen-decr}, \req{wik-bound-ming} and \req{viklow} that
\begin{align}
f(x_0)-f(x_{j_\eta+1})
&\ge \sum_{j=0}^{j_\eta} \sum_{i=1}^n \frac{\tau\varsigma_{\min}g_{i,j}^2}{2\kB w_{i,j}}
- \half (\kB+L)\sum_{j=0}^{j_\eta} \sum_{i=1}^n\frac{g_{i,j}^2}{w_{i,j}^2}\neol
&\ge - \half (\kB+L)\sum_{j=0}^{j_\eta}\sum_{i=1}^n\frac{g_{i,j}^2}{w_{i,j}^2}\neol
&\ge - \half (\kB+L)\sum_{j=0}^{j_\eta}\sum_{i=1}^n\frac{g_{i,j}^2}{\max[\varsigma,v_{i,k}]^2(j+1)^{2\nu} \theta^2}\neol
&\ge - \half \frac{(\kB+L)}{\theta^2}\sum_{j=0}^{j_\eta}\sum_{i=1}^n\frac{g_{i,j}^2}{v_{i,k}^2(j+1)^{2\nu}}\neol
&\ge - \half  \frac{(\kB+L)n}{\theta^2}\sum_{j=0}^{j_\eta} h(j)^2.\label{A2}
\end{align}
Combining \req{A1} and \req{A2}, we obtain that
\[
\Gamma_0
\ge f(x_0)-f(x_{k+1})
\ge - \half \frac{(\kB+L)n}{\theta^2}\sum_{j=0}^{j_\eta} h(j)^2 + \frac{\eta (k-j_\eta)}{2\kB(k+1)^\mu \theta}
     \min_{j\in\iibe{j_\theta+1}{k}}\left(\sum_{i=1}^n\min\left[\frac{g_{i,j}^2}{\varsigma_i},\frac{g_{i,j}^2}{v_{i,j}}\right]\right)
\]
and thus that
\[
\min_{j\in\iibe{j_\theta+1}{k}}\left(\sum_{i=1}^n\min\left[\frac{g_{i,j}^2}{\varsigma_i},\frac{g_{i,j}^2}{v_{i,j}}\right]\right)
\le \frac{2\kB(k+1)^\mu}{\eta (k-j_\eta)}\left[\Gamma_0 \theta +\half \frac{n (\kB+L)}{\theta}\sum_{j=0}^{j_\eta} h(j)^2\right]
\]
and we deduce that there must exist a subsequence
$\{k_\ell\}\subseteq \{k\}_{j_\eta+1}^\infty$ such that, for each $\ell$,
\beqn{A3}
\sum_{i=1}^n\min\left[\frac{g_{i,k_\ell}^2}{\varsigma_i},\frac{g_{i,jk_\ell}^2}{v_{i,k_\ell}}\right]
\le \frac{2\kB(k_\ell+1)^\mu}{\eta (k_\ell-j_\eta)}\left[\Gamma_0 \theta +\half \frac{n (\kB+L)}{\theta}\sum_{j=0}^{j_\eta} h(j)^2\right].
\eeqn
But
\beqn{A0}
\frac{(k_\ell+1)^\mu}{k_\ell-j_\eta}
< \frac{2^\mu k_\ell^\mu}{k_\ell-j_\eta}
< \frac{2 k_\ell^\mu}{k_\ell-j_\eta}
=\frac{2 k_\ell^\mu k_\ell}{(k_\ell-j_\eta)k_\ell}
=\frac{k_\ell}{k_\ell-j_\eta} \cdot \frac{2}{k_\ell^{1-\mu}}
\le  \frac{2(j_\eta+1)}{k_\ell^{1-\mu}},
\eeqn
where we used the facts that $\mu <1$ and that $\frac{k_\ell}{k_\ell-j_\theta} $
is a decreasing function for $k_\ell \ge j_\theta+1$.
Using this inequality, we thus obtain from \req{A3} that, for each $\ell$,
\[
\sum_{i=1}^n\min\left[\frac{g_{i,k_\ell}^2}{\varsigma_i},\frac{g_{i,jk_\ell}^2}{v_{i,k_\ell}}\right]
\le \frac{4\kB(j_\eta+1)}{\eta \, k_\ell^{1-\mu}}\left[\Gamma_0 \theta
  +\half n \frac{(\kB+L)}{\theta}\sum_{j=0}^{j_\eta} h(j)^2\right].
\]
As a consequence,
\[
k_\varsigma \eqdef
\left(\frac{4\kB(j_\eta+1)\left[\Gamma_0 \theta
      +\half \frac{n (\kB+L)}{\theta}\sum_{j=0}^{j_\eta} h(j)^2\right]}{\eta\varsigma_{\min}}\right)^\frac{1}{1-\mu}
\]
is such that, for all $k_\ell \ge k_\varsigma$,
\[
\min\left[\frac{g_{i,k_\ell}^2}{\varsigma_i,},\frac{g_{i,_\ell}^2}{v_{i,k\ell}}\right]
\le  \varsigma_{\min}.
\]
Lemma~\ref{lemma:divsprop} then yields that, for all $k_\ell \ge k_\varsigma$,
\begin{align*}
\sum_{i=1}^n \frac{g_{i,k_\ell}^2}{2 \varsigma_i}
&\le \frac{2\kB(k_\ell+1)^\mu}{\eta (k_\ell-j_\eta)}\left[\Gamma_0 \theta +\half \frac{n (\kB+L)}{\theta}\sum_{j=0}^{j_\eta} h(j)^2\right]
\end{align*}
which, because $\varsigma_i \le 1$, gives that, for all $k_\ell \ge k_\varsigma$,
\beqn{choosetheta}
\|g_{k_\ell}\|^2
\le
\frac{(k_\ell+1)^\mu}{k_\ell-j_\eta}\,
\left(\frac{4\kB}{\eta}\right)
  \left[\Gamma_0 \theta +\half n \frac{(\kB+L)}{\theta}\sum_{j=0}^{j_\eta} h(j)^2\right],
\eeqn
finally implying \req{ngkbound-ming} because of \req{A0}.
} 

\noindent
We again provide some comments on this last result.
\begin{enumerate}
\item
  The choice \req{wik-bound-ming} is of course reminiscent, in a smooth and
  nonconvex setting, of the ``diminishing stepsize'' subgradient method for stochastic problems (see \cite[Theorem~1.2.4]{Bert95} or \cite[Theorems~8.25 and 8.40]{Beck17} 
  and the many references therein), for which a $\calO(1/\sqrt{k})$ global rate of convergence is
  known.
\item
  Theorem~\ref{theorem:ming-style} provides information on the speed
  of convergence for iterations that are beyond an \emph{a priori}
  computable iteration index. Indeed $j_\theta$ and $k_\varsigma$ only
  depends on $\nu$ $h(\ell)$  and problem's constants and, in
  particular do not depend on $k$. However, the formulation of the
  theorem is slightly weaker than that of
  Theorem~\ref{theorem:allmu}. Because \req{ngkbound-ming} only holds
  for iterates along the subsequence $\{k_\ell\}$, there is no
  guarantee that the bounnd given by the right-hand-side is valid at
  other iterations. But note that this right-hand side depends on
  $k_\ell$, which is an index in the complete sequence of iterates,
  rather than on $\ell$ (the subsequence index).
  
  This slightly weaker formulation is no longer necessary if one is
  ready to assume bounded gradients, as can be seen in Theorem~4.1 in \cite{GratJeraToin22b}.  
\item  As the chosen values of $\mu$ and $\nu$ approach zero, then
  the $k$-order of convergence
  beyond $j_\theta$ tends to $\calO(1/\sqrt{k_\ell})$,
  which the order derived for the methods of the previous section and is the standard $k$-order
  for first-order methods using evaluations of the objective function, albeit the value
  of $j_\theta$ might increase. 
\item As in Section~\ref{adap-s} and considering \eqref{choosetheta}, we may choose $\theta$ according to \req{thetabest}
in an attempt to balance the two terms of the letf-hand side and improve the explicit dependence of 
the complexity bound on $n$. 
\end{enumerate}

\noindent
We are now again interested to estimate how sharp the $k$-order bound
\req{ngkbound-ming} in $\calO(\frac{1}{k^{(1-\mu)/2}})$ is.

\lthm{sharp2}{
The bound \req{ngkbound-ming} is essentially sharp in that, for any $\omega >
\half(1-\nu)$, there exists a univariate function $f_\omega(x)$ satisfying AS.1--AS.3
such that the \al{ASTR1} algorithm with \req{wik-bound-ming}, $B_k=0$ and $\theta=1$ applied to this
function produces a sequence of gradient norms given by $\|g_k\| = \frac{1}{(k+1)^\omega}$.
}

\proof{
Consider the sequence defined, for some $\omega \in(\half(1-\nu),1]$  and all
$k\geq 0$, by 
\beqn{seqs}
g_k = - \frac{1}{(k+1)^\omega}
\ms
w_k
= \max\left[\varsigma,\max_{\ell\in\iiz{k}}|g_\ell|\right]\,(k+1)^{\nu}
= (k+1)^{\nu},
\eeqn
\beqn{seqsb}
s_k = \frac{1}{(k+1)^{2\omega-\nu}}<1
\tim{ and }
f_{k+1} = f_k+g_ks_k,
\eeqn
where we have chosen $\varsigma\in (0,1)$ and $f_0 = \zeta(2\omega+\half)$ where $\zeta(\cdot)$ is the
Riemann zeta function. Immediately note that 
\[
\lim_{k\rightarrow\infty}|g_k| = 0,
\]
and $|g_k|\leq=1=\kappa_g$ for all $k$. We now verify that, if
\[
x_0 = 0 \tim{ and } x_k = x_{k-1} + s_{k-1} \tim{for} k \geq 1,
\]
then exists a function $f_\omega(x)$ satisfying AS.1--AS.3 such that, for all $k\ge0$,
\[
f_\omega(x_k) = f_k,
\tim{ and }
g_\omega(x_k) = g_k,
\]
and such that the sequence defined by \req{seqs}-\req{seqsb} is generated by applying the
\al{ASTR1} algorithm using $B_k=0$ and $\theta=1$.  The function $f_\omega(x)$ is constructed using Hermite
interpolation on each interval $[x_k,x_{k+1}]$ (note that the $x_k$ are
monotonically increasing), which known (see \cite{CartGoulToin10a} or
\cite[Th. A.9.2]{CartGoulToin22}) to
exist whenever there exists a constant $\kappa_f\ge 0$ such that, for each $k$,
\[
|f_{k+1}-f_k-g_ks_k| \le \kappa_f |s_k|^2
\tim{ and }
|g_{k+1}- g_k | \le \kappa_f |s_k|.
\]
The first of these conditions holds by construction of the
$\{f_k\}_{k\ge0}$. To verify the second, we first note that,
because $1/(k+1)^\omega$ is a convex function of $k$ and $|1/(k+1)|\le 1$,
\beqn{exg1}
\frac{|g_{k+1}-g_k|}{|s_k|}
\leq \frac{\omega (k+1)^{2\omega-\nu}}{(k+1)^{1+\omega}}
= \frac{\omega}{(k+1)^{\nu-\omega+1}}
\leq \omega
\ms
(k \geq 0),
\eeqn
where $\nu-\omega+1 \ge \nu >0$, so that the desired inequality holds with
$\kappa_f=\omega$.

Moreover, Hermite interpolation guarantees that $f_\omega(x)$ is bounded below
whenever $|f_k|$ and $|s_k|$ remain bounded. We have already verified the
second of these conditions in \req{seqsb}. We also have from \req{seqsb} that
\beqn{df-sh2}
f_0-f_{k+1}
= \sum_{j=0}^k \frac{1}{j+1)^{2\omega}(j+1)^{\nu}}
\eeqn
which converges to the finite limit $\zeta(2\omega+\nu)$ because we have
chosen $\omega> \half(1-\nu)$.  Thus $f_k \in (0,\zeta(2\omega+\nu)]$ for all
$k$ and the first condition is also satisfied and AS.3 holds.
This completes our proof.
}

\noindent
The conclusions which can be drawn from this theorem parallel those drawn
after Theorem~\ref{sharp1}. The bound \req{ngkbound-ming} is essentially sharp
(in the sense of \cite{CartGoulToin18a}\footnote{Observe that $f_0$ now tends
to infinity when $\omega$ tends to $\half(\nu-1)$ and hence that AS.3 fails in the
limit. As before, the structure of \req{ngkbound-ming} implies that the
complexity bound deteriorates when the gap $\Gamma_0 = f(x_0)-\flow$ grows.}) for
the \al{ASTR1} algorithm with \req{wik-bound-ming}.


\numsection{Numerical illustration}\label{numerics-s}

We now provide some numerical illustration on problems which are commonly used for 
the evaluation of optimization algorithms. For the sake of clarity and conciseness,
we needed to keep the list of algorithmic variants reported here reasonably limited,
and have taken the following considerations into account for our choice.
\begin{enumerate}
\item Both weights' definitions \req{w-adag} and \req{wik-bound-ming} are illustrated.
  Moreover, since the Adam algorithm using \req{w-adam} is so commonly used the
  stochastic context, we also included it in the comparison.
\item Despite Theorems~\ref{theorem:allmu} and \ref{theorem:ming-style} covering a
  wide choice of the parameters $\mu$ and $\nu$, we have chosen to focus here on the
  most common choice for \req{w-adag} and \req{w-adam} (i.e. $\mu = \half$ and
  $\beta_2=\sfrac{9}{10}$, corresponding to Adagrad, Adagrad-Norm and Adam). When 
  using \req{wik-bound-ming}, we have also restricted our comparison to the single 
  choice of $\mu$ and $\nu$ used (with reasonable success) in \cite{GratJeraToin22a}, 
  namely $\mu=\nu=\sfrac{1}{10}$.
   By contrast, we have included results for the variants of the weights' definitions
  for values of $\theta=1$ (the standard choice) and $\theta = \sqrt{n}$  suggested in 
  \req{thetabest} for unknown $\Gamma_0$.
\item In order to be able to test enough algorithmic variants on enough problems in
  reasonable computing time, we have decided to focus our experiments on
  low-dimensional problems in the case where $\vartheta=1$. We have nevertheless considered a few 
  large-scale instance for the purpose of illustrating the effect of the scaling 
  $\theta = \sqrt{n}$ which is only expected to be significant for such instances.
\item We have chosen to define the step $s_k$ in Step~3 of the \al{ASTR1} algorithm
  by approximately minimizing the quadratic model \req{q-model} within the
  $\ell_\infty$ trust-region using a projected truncated conjugate-gradient
  approach \cite{MoreTora89,MoreTora91} which is terminated as soon as 
  \[
  \|g_k+B_ks_k\|_2\leq \max\Big[10^{-12},10^{-5} \|g_k\|_2\Big].
  \]
  We also considered an alternative,
  namely that of minimizing the quadratic model in an Euclidean $\ell_2$ trust region
  (with the same accuracy requirement) using a Generalized Lanczos Trust
  Region (GLTR) technique \cite{GoulLuciRomaToin99}.
\item We thought it would be interesting to compare ``purely first-order''
  variants (that is variants for which $B_k=0$ for all $k$) with methods
  using some kind of Hessian approximation. Among many possibilities, we selected
  three types of approximations of interest. The first is the diagonal Barzilai-Borwein
  approximation \cite{BarzBorw88} 
  \beqn{BBH}
  B_{k+1} = \frac{\|s_k\|_2^2}{y_k^Ts_k} I_n
  \eeqn
  where $I_n$ is the identity matrix of dimension $n$, $y_k=g_{k+1}-g_k$ and
  $y_k^Ts_k \geq 10^{-15}\|s_k\|_2^2$.  The second 
  is limited-memory BFGS approximations \cite{LiuNoce89}, where a small number (3)
  of BFGS updates are added to the matrix \req{BBH}, each update corresponding
  to a secant pair $(y_k,s_k)$ with $y_k^Ts_k \geq 10^{-15}\|s_k\|_2^2$. The third is not to approximate the Hessian
  at all, but to use its exact value, that is $B_k=\nabla_x^2f(x_k)$ for all $k$.
\end{enumerate}
Given these considerations, we have selected the algorithmic
\al{ASTR1} variants using $\gamma_k=1$ and detailed in
Table~\ref{variants}, where the second column indicates the norm used
to define the trust-region.  
 
\begin{table}[t]
\centerline{
{\renewcommand{\arraystretch}{1.7}%
\begin{tabular}{|l|l|l|c|l|}
\hline
Name          & Norm               & \multicolumn{1}{c|}{$w_{i,k}$ definition ($i\in\ii{n}$)} & $B_{k+1}$ & params \\
\hline
{\tt adagnorm}   & $\|\cdot\|_2$      & $w_{i,k} = \left[\sfrac{1}{100}+\sum_{j=0}^k\|g_j\|_2^2\right]^\half$ & 0 & \\
\hline
{\tt adagrad} & $\|\cdot\|_\infty$ & $w_{i,k} = \left[\sfrac{1}{100}+\sum_{j=0}^k g_{i,j}^2\right]^\half$ & 0 & \\
\hline
{\tt adamnorm}   & $\|\cdot\|_2$      & $w_{i,k} = \left[\sfrac{1}{100}+\sum_{j=0}^k\beta_2^{k-j}\|g_j\|_2^2\right]^\half$ & 0 & $\beta = \sfrac{9}{10}$ \\
\hline
{\tt adam}    & $\|\cdot\|_\infty$ & $w_{i,k} = \left[\sfrac{1}{100}+\sum_{j=0}^k\beta_2^{k-j}g_{i,j}^2\right]^\half$ & 0 & $\beta = \sfrac{9}{10}$ \\
\hline
{\tt maxgnorm}  & $\|\cdot\|_2$      & $w_{i,k} = (k+1)^\sfrac{1}{10}\max\Big[ \sfrac{1}{100},\bigmax_{j\in\iiz{k}}\|g_j\|_2 \Big]$ & 0 & \\
\hline
{\tt maxg} & $\|\cdot\|_\infty$ & $w_{i,k} = (k+1)^\sfrac{1}{10}\max\Big[ \sfrac{1}{100},\bigmax_{j\in\iiz{k}}|g_{i,j}| \Big]$ & 0 &\\
\hline
{\tt adagbb} & $\|\cdot\|_\infty$ & $w_{i,k} = \left[\sfrac{1}{100}+\sum_{j=0}^k g_{i,j}^2\right]^\half$ & \req{BBH} & \\
\hline
{\tt adagbfgs3} &$\|\cdot\|_\infty$ & $w_{i,k} = \left[\sfrac{1}{100}+\sum_{j=0}^k g_{i,j}^2\right]^\half$ & LBFGS(3) & \\
\hline
{\tt adagH} & $\|\cdot\|_\infty$ & $w_{i,k} = \left[\sfrac{1}{100}+\sum_{j=0}^k g_{i,j}^2\right]^\half$ & $\nabla_x^2f(x_{k+1})$ & \\
\hline
{\tt adagrads} & $\|\cdot\|_\infty$ & $w_{i,k} = \sqrt{n}\left[\sfrac{1}{100}+\sum_{j=0}^k g_{i,j}^2\right]^\half$ & 0 & \\
\hline
{\tt adams} & $\|\cdot\|_\infty$ & $w_{i,k} = \sqrt{n}\left[\sfrac{1}{100}+\sum_{j=0}^k\beta_2^{k-j}g_{i,j}^2\right]^\half$ & 0 & $\beta = \sfrac{9}{10}$ \\
\hline
{\tt maxgs} & $\|\cdot\|_\infty$ & $w_{i,k} = \sqrt{n}(k+1)^\sfrac{1}{10}\max\Big[ \sfrac{1}{100},\bigmax_{j\in\iiz{k}}|g_{i,j}| \Big]$ & 0 &\\
\hline
{\tt adagbbs} & $\|\cdot\|_\infty$ & $w_{i,k} = \sqrt{n}\left[\sfrac{1}{100}+\sum_{j=0}^k g_{i,j}^2\right]^\half$ & \req{BBH} & \\
\hline
{\tt adagbfgs3s} &$\|\cdot\|_\infty$ & $w_{i,k} = \sqrt{n}\left[\sfrac{1}{100}+\sum_{j=0}^k g_{i,j}^2\right]^\half$ & LBFGS(3) &\\
\hline
{\tt adagHs} & $\|\cdot\|_\infty$ & $w_{i,k} = \sqrt{n}\left[\sfrac{1}{100}+\sum_{j=0}^k g_{i,j}^2\right]^\half$ & $\nabla_x^2f(x_{k+1})$ & \\
\hline
{\tt sdba}    & \multicolumn{4}{l|}{standard steepest-descent with
  backtracking (e.g. \cite[Algorithm~2.2.1]{CartGoulToin22})} \\
\hline
\end{tabular}}
}
\caption{\label{variants}The considered algorithmic variants}
\end{table}

Note that all variants for which $B_{k+1} = 0$ are "purely first-order" in the sense discussed above.
Note also that, under AS.3, {\tt maxgnorm} and {\tt maxg} satisfy \req{wik-bound-ming}
with $\mu=\nu=\sfrac{1}{10}$, $\varsigma_i=\varsigma=\sfrac{1}{100}$ and $\kappa_w=\kappa_g$. 
All algorithms were run\footnote{In
Matlab\textregistered\ running under Ubuntu on a Dell Precision with 16 cores and 64
GB of memory.} on the low dimensional instances of the
problems\footnote{From their standard starting point.} of the {\sf{OPM}}
collection \cite{GratToin21c} (April 2023), a subset of widely used CUTEst testing environment \cite{GoulOrbaToin15b}. 
The instances are listed with their dimension in
Table~\ref{testprobs}, until either $\|\nabla_x^1f(x_k)\|_2\leq 10^{-6}$, or a
maximum of 100000 iterations was reached, or evaluation of the derivatives
returned an error.

\begin{table}[htb]\footnotesize
\begin{center}
  \begin{tabular}{|l|r|l|r|l|r|l|r|l|r|l|r|}
    \hline
Problem & $n$ & Problem & $n$ & Problem & $n$ & Problem & $n$ & Problem & $n$ & Problem & $n$ \\
\hline
argauss       &  3 & chebyqad    & 10 & dixmaanl    & 12 & heart8ls   &  8 & msqrtals    & 16 & scosine     & 10 \\ 
arglina       & 10 & cliff       &  2 & dixon       & 10 & helix      &  3 & msqrtbls    & 16 & sisser      &  2 \\
arglinb       & 10 & clplatea    & 16 & dqartic     & 10 & hilbert    & 10 & morebv      & 12 & spmsqrt     & 10 \\
arglinc       & 10 & clplateb    & 16 & edensch     & 10 & himln3     &  2 & nlminsurf   & 16 & tcontact    & 49 \\
argtrig       & 10 & clustr      &  2 & eg2         & 10 & himm25     &  2 & nondquar    & 10 & trigger     &  7 \\
arwhead       & 10 & cosine      & 10 & eg2s        & 10 & himm27     &  2 & nzf1        & 13 & tridia      & 10 \\
bard          &  3 & crglvy      &  4 & eigfenals   & 12 & himm28     &  2 & osbornea    &  5 & tlminsurfx  & 16 \\
bdarwhd       & 10 & cube        &  2 & eigenbls    & 12 & himm29     &  2 & osborneb    & 11 & tnlminsurfx & 16 \\
beale         &  2 & curly10     & 10 & eigencls    & 12 & himm30     &  3 & penalty1    & 10 & vardim      & 10 \\
biggs5        &  5 & dixmaana    & 12 & engval1     & 10 & himm32     &  4 & penalty2    & 10 & vibrbeam    &  8 \\
biggs6        &  6 & dixmaanb    & 12 & engval2     &  3 & himm33     &  2 & penalty3    & 10 & watson      & 12 \\
brownden      &  4 & dixmaanc    & 12 & expfit      &  2 & hypcir     &  2 & powellbs    &  2 & wmsqrtals   & 16 \\
booth         &  2 & dixmaand    & 12 & extrosnb    & 10 & indef      & 10 & powellsg    & 12 & wmsqrtbls   & 16 \\
box3          &  3 & dixmaane    & 12 & fminsurf    & 16 & integreq   & 10 & powellsq    &  2 & woods       & 12 \\
brkmcc        &  2 & dixmaanf    & 12 & freuroth    &  4 & jensmp     &  2 & powr        & 10 & yfitu       &  3 \\
brownal       & 10 & dixmaang    & 12 & genhumps    &  5 & kowosb     &  4 & recipe      &  2 & zangwill2   &  2 \\
brownbs       &  2 & dixmaanh    & 12 & gottfr      &  2 & lminsurg   & 16 & rosenbr     & 10 & zangwill3   &  3 \\
broyden3d     & 10 & dixmaani    & 12 & gulf        &  4 & macino     & 10 & sensors     & 10 &             &    \\
broydenbd     & 10 & dixmaanj    & 12 & hairy       &  2 & mexhat     &  2 & schmvett    &  3 &             &    \\
chandheu      & 10 & dixmaank    & 12 & heart6ls    &  6 & meyer3     &  3 & scurly10    & 10 &             &    \\
\hline
\end{tabular}
\caption{\label{testprobs} The small \sf{OPM} test problems and their dimension}
\end{center}
\end{table}

Before considering the results, we make two additional comments.  The first is
that very few of the test functions have bounded gradients on the whole of $\Re^n$.
While this is usually not a problem when testing standard first-order descent
methods (because it may then be true in the level set determined by the
starting point), this is no longer the case for significantly non-monotone
methods like the ones tested here.  As a consequence, it may (and does) happen
that the gradient evaluation is attempted at a point where its value exceeds
the Matlab overflow limit, causing the algorithm to fail on the
problem.  The second comment is that the (sometimes quite wild)
non-monotonicity of the methods considered here has another practical
consequence: it happens on several nonconvex problems\footnote{broyden3d,
broydenbd, curly10, gottfr, hairy, indef, jensmp, osborneb, sensors,
wmsqrtals, wmsqrtbls, woods. } that convergence of different algorithmic
variants occurs to points with gradient norm within termination tolerance (the
methods are thus achieving their objective), but these points can be quite far
apart and may have very different function values.  

We discuss the results of our tests from the efficiency and reliability points
of view. Efficiency is measured in number of derivatives' evaluations
(or, equivalently, iterations)\footnote{For {\tt sdba}, gradient and
objective-function evaluations.}: the fewer evaluations the more efficient the
algorithm. Because the standard performance profiles
\cite{DolaMoreMuns06} for our selection of 16 algorithms would be too crowded to read, 
we follow \cite{PorcToin17c} and consider the derived
``global'' measure $\pi_{\tt algo}$ to be $\sfrac{1}{50}$ of the area below the curve
corresponding to {\tt algo} in this performance profile, for abscissas in the
interval $[1,50]$. The larger this area and closer $\pi_{\tt algo}$ to one,
the closer the curve to the right and top borders of the plot and the better
the global performance. When reporting reliability, we say that the run of an 
algorithmic variant on a specific test problem is successful if the gradient
norm tolerance has been achieved.
In what follows, $\rho_{\tt algo}$ denotes the percentage of successful runs taken on
all problems.  Table~\ref{nonoise} presents
the values of these statistics in two columns: for easier reading, the
variants are sorted by decreasing global performance ($\pi_{\tt algo}$) in the first, and
by decreasing reliability ($\rho_{\tt algo}$) in the second.
\noindent
A total of 18 problems\footnote{
  biggs5,         brownbs,        cliff,          genhumps,       gulf,
  heart6ls,       heart8ls,       himm29,         mexhat,         meyer3,
  nondquar,       osbornea,       penalty2,       powellbs,       powellsg,
  scurly10,       watson,         yfitu.          
} could not be successfully solved by any of the above algorithms, we believe
mostly because of ill-conditioning.

\comment{
\begin{figure}[htb] 
\centerline{
\includegraphics[width=10cm]{./gjt2R4-noiseless.eps}
}
\caption{\label{fig:profile}Performance profile for deterministic OFFO algorithms on
  {\sf{OPM}} problems}
\end{figure}
}
\begin{table}
\begin{center}
\begin{tabular}{|l|r|r||l|r|r|}
\hline
Method & $\pi_{\tt algo}$  & $\rho_{\tt algo}$ & {\tt algo} & $\pi_{\tt algo}$ & $\rho_{\tt algo}$   \\
\hline
{\tt adagbfgs3}   &  0.75  &  69.75   &  {\tt adagrad}    &   0.69  &   73.11 \\
{\tt sdba}        &  0.73  &  68.91   &  {\tt adagbfgs3}  &   0.75  &   69.75 \\
{\tt adagH}       &  0.72  &  69.75   &  {\tt adagH}      &   0.72  &   69.75 \\
{\tt adagrad}     &  0.69  &  73.11   &  {\tt sdba}       &   0.73  &   68.91 \\
{\tt maxg}        &  0.66  &  66.39   &  {\tt adagHs}     &   0.63  &   67.23 \\
{\tt adagHs}      &  0.63  &  67.23   &  {\tt maxg}       &   0.66  &   66.39 \\
{\tt adagbb}      &  0.63  &  64.71   &  {\tt adagrads}   &   0.60  &   65.55 \\
{\tt adagbfgs3s}  &  0.62  &  63.87   &  {\tt adagbb}     &   0.63  &   64.72 \\
{\tt maxgs}       &  0.60  &  62.18   &  {\tt adagbfgs3s} &   0.62  &   63.87 \\
{\tt adagrads}    &  0.59  &  65.55   &  {\tt maxgs}      &   0.60  &   62.18 \\
{\tt adagnorm}    &  0.58  &  61.34   &  {\tt adagnorm}   &   0.58  &   61.34 \\
{\tt maxgnorm}    &  0.56  &  57.98   &  {\tt adagbbs}    &   0.56  &   60.50 \\
{\tt adagbbs}     &  0.56  &  60.50   &  {\tt maxgnorm}   &   0.56  &   57.98 \\
{\tt adamnorm}    &  0.55  &  34.45   &  {\tt adamnorm}   &   0.55  &   34.45 \\
{\tt adam}        &  0.54  &  30.25   &  {\tt adams}      &   0.52  &   33.61 \\
{\tt adams}       &  0.52  &  33.61   &  {\tt adam}       &   0.54  &   30.25 \\
\hline
\end{tabular}
\caption{\label{nonoise} Performance and reliability statistics for deterministic OFFO and steepest descent algorithms on
small {\sf{OPM}} problems ($\epsilon_1=10^{-6}$)}
\end{center}
\end{table}

The authors are of course aware that the very limited experiments presented
here do not replace extended numerical practice and could be completed in
various ways. They nevertheless suggest the following (very tentative)
comments.

\begin{enumerate}
\item There often seems to be a definite advantage in using the $\|\cdot\|_\infty$
  norm over $\|\cdot\|_2$, as can be seen by comparing {\tt adagnorm} with {\tt
    adagrad} and {\tt maxgnorm} with {\tt
    maxg}. While this may be due in part to the fact that the trust region
  in $\ell_\infty$ norm is larger than that in $\ell_2$ norm (and thus allows
  larger steps), it is also the case that the disaggregate definition of the
  scaling factors $w_{i,k}$ (\req{w-adag}, \req{w-adam} or
  \req{wik-bound-ming}) used in conjunction with the $\ell_\infty$ norm may
  allow a better exploitation of differences of scale between coordinates.
\item Among the ''purely first-order'' methods, {\tt sdba}, {\tt maxg} and
  {\tt adagrad} are almost undistinguishable form the performance point of
  view, with a reliability advantage for {\tt adagrad} (the most reliable 
  method in our tests). This means that, at least in those experiments, the 
  suggestion resulting from the theory that OFFO methods may perform 
  comparably to standard first-order methods seems vindicated.
\item The Adam variants ({\tt adamnorm} and {\tt adam}) are clearly
  outperformed in our tests by the Adagrad ones ({\tt adagnorm} and {\tt
  adagrad}). We recall that analytical examples where Adam fails do exist,
  while the convergence of Adagrad is guaranteed.
\item The theoretical difference in global rate of convergence between {\tt
  adagrad} and {\tt maxg} does not seem to have much impact on the relative
  performance of these two methods.
\item The use of limited memory Hessian approximation ({\tt adagbfgs3}) appears
  to enhance the performance of {\tt adagrad}, but this is not the case of the
  Barzilai-Borwein approximation ({\tt adagbb}) or, remarkably, for the use
  of the exact Hessian ({\tt adagH}).  When these methods fail, this is
  often  because the steplength is too small to allow the truncated
  conjugate-gradient solver to pick up second-order information in other 
  directions than the negative gradient. What favours the limited memory
  approach remains unclear at this stage.
\end{enumerate}

We also note that the variants scaled with \req{thetabest} (with
$\Gamma_0$ unknown) did not perform better on small dimensional
problems; possibly because the factor $n$ does not dominate the
complexity bounds in this case.  To illustrate 
the impact of this scaling in larger cases, we ran a subset of six
methods which performed well in Table~\ref{nonoise} (namely {\tt
  adagnorm}, {\tt adagrad}, {\tt adagrads}, {\tt maxgnorm}, {\tt maxg}
and {\tt maxgs}) on the broyden3d and nlminsurf problems with
increasing problem dimension (we used $\epsilon_1=10^{-3}$).  
The results are reported in Table~\ref{table:large}.

\begin{table}[h!]\small
\begin{tabular}{|lr|r|r|r|r|r|r|r|r|r|r|}
\hline
Method  &     & \multicolumn{5}{c|}{broyden3d}&\multicolumn{5}{c|}{nlminsurf}\\
        & $n$ &       10  &   100 &  1000 &  10000 & 100000 &   256 &  1034 &   4096 &  16384 &  65536 \\
\hline
{\tt adagnorm}  &&     37   &    71 &   467 &   4257 &  43400 &    166 &   503 &   1791 &   6038 &  19239 \\
{\tt adagrad}   &&    200   & 37809 & 37809 &  37809 &  37809 &   7966 & 30795 & 121164 & 482025 &     NR \\
{\tt adagrads}  &&    134   &   190 &  1452 &  13042 & 125503 & 138 &   532 &   2981 &  19414 & 121934 \\
{\tt maxgnorm}  &&     46   &    76 &   285 &   1138 &   4520 & 1699 &  3978 &   3867 &   5355 &  19424 \\
{\tt maxg}      &&    458   &   410 &   462 &   3362 &  36609 &   NC &    NC &     NC &     NC &     NR \\
{\tt maxgs}     &&     76   &   155  &  567 &   2048 &   7370 &  1142 &  1155 &   5049 &   6407 &  30661 \\
\hline
\end{tabular}
\caption{Number of iterations for convergence on the broyden3d and
  nlminsurf problems as a function of dimension ($\epsilon_1=10^{-3}$,
  NC = more than $10^6$ iterations, NR = not run)\label{table:large}} 
\end{table}

Despite the improvement in complexity due to choosing $\theta>1$ being
theoretical (and applies to the worst-case performance), we may still
note some positive (if not completely uniform) effect in this table.
The interpretation is also blurred somewhat by the fact that {\tt
  maxg} and {\tt adagrad} converged to local minimas of broyden3d
rather than the global one. We nevertheless note the consistently
better performance of {\tt adagnorm} compared to {\tt adagrad} and
{\tt adagrads}, possibly illustrating the fact that its complexity
bound does not explicitly involve $n$. 

Finally, and although this is a slight digression from the paper's main topic,
we report in Table~\ref{noisy} how reliability of our selection of OFFO
variants is impacted by noise.  To obtain these results, we ran the considered
methods on all test problems where the evaluations (function\footnote{For {\tt
    sdba}.} and derivatives) are contaminated by 5, 15, 25 or 50 \%\ of
relative Gaussian noise with unit variance.  The reliability percentages in
the table result from averaging results obtained for ten independent runs. 

\begin{table}[htb]
\begin{center}
\begin{tabular}{|l|r|r|r|r|r|}
\hline
                 &\multicolumn{5}{c|}{$\rho_{\tt algo}$/relative noise level}\\
\cline{2-6}
{\tt algo}       &  0\%  &   5\%  &  15\% & 25\% & 50\% \\
\hline
{\tt adagH}      & 83.19 & 84.96 & 84.20 & 84.71 & 82.18 \\
{\tt adagHs}     & 81.51 & 81.85 & 81.91 & 80.50 & 77.82 \\
{\tt adagbfgs3}  & 78.15 & 80.50 & 80.50 & 80.84 & 80.18 \\
{\tt adagrad}    & 77.31 & 80.50 & 80.25 & 80.17 & 80.17 \\
{\tt adagbb}     & 75.69 & 80.08 & 80.17 & 79.58 & 79.41 \\
{\tt adagbfgs3s} & 78.99 & 79.50 & 70.67 & 79.41 & 78.66 \\
{\tt adagbbs}    & 73.95 & 78.15 & 78.40 & 78.49 & 77.06 \\
{\tt adagrads}   & 78.15 & 78.07 & 78.66 & 78.74 & 77.23 \\
{\tt maxgs}      & 75.63 & 76.39 & 75.46 & 76.05 & 74.54 \\
{\tt adagnorm}   & 75.63 & 75.21 & 75.80 & 75.71 & 74.03 \\
{\tt maxg}       & 74.79 & 74.37 & 75.55 & 78.15 & 78.07 \\
{\tt maxgnorm}   & 69.75 & 68.74 & 69.75 & 70.84 & 71.01\\
{\tt adams}      & 42.86 & 37.98 & 40.25 & 44.79 & 50.84 \\
{\tt adamnorm}   & 42.02 & 37.56 & 44.96 & 50.84 & 55.29 \\
{\tt adam}       & 40.34 & 35.55 & 36.30 & 44.03 & 45.80 \\
{\tt sdba}       & 81.51 & 30.92 & 31.85 & 34.87 & 29.58 \\
\hline
\end{tabular}
\caption{\label{noisy} Reliability of OFFO algorithms and steepest descent as a function of
  the level of relative Gaussian noise ($\epsilon_1 = 10^{-3}$)}
\end{center}
\end{table}

As can be seen in the table, the reliability of the {\tt sdba} methods
dramatically drops as soon as noise is present, while that of the other OFFO
methods is barely affected and remains globally unchanged\footnote{
It is interesting that reliability is slightly better for the noisy 
cases and the better OFFO methods.} 
for increasing noise
level.  This is consistent with widespread experience in the deep learning
context, where noise is caused by sampling among the very large number of
terms defining the objective function.  This observation vindicates
the popularity of methods such as Adagrad in the noisy context and suggests
that the new OFFO algorithms may have some practical potential.

We conclude by noting that the algorithms' reliability ($\rho_{\tt algo}$ 
is (expectedly) better for $\epsilon_1 = 10^{-3}$ (first column of 
Table~\ref{noisy}) than for $\epsilon_1 = 10^{-6}$  (Table~\ref{nonoise}), 
but that the improvement remains modest, the reliability of Adagrad
decreasing marginally slower. 

\numsection{Conclusions}\label{concl-s}

We have presented a parametric class of deterministic ``trust-region minded''
extensions of the Adagrad method, allowing for the use of second-order
information, should it be available.  We then prove that, for OFFO algorithms
in this class, $\min_{j\in\iiz{k}}\|g_j\| = \calO(1/\sqrt{k+1})$.
We have also shown that this bound, which does not require any uniform bound on the gradient,
is essentially sharp. It is \emph{identical to the global
rate of convergence of standard first-order methods using both
objective-function and gradient evaluations}, despite the fact that the
latter exploit significantly more information. Thus, \emph{if one considers
the order of global convergence only, evaluating the objective-function values
is an unnecessary effort.}  We have also considered another class of OFFO
algorithms inspired by the ``diminishing stepsize'' paradigm in non-smooth
convex optimization and have provided an essentially sharp (but slighlty
worse) global rate of convergence for this latter class. Limited numerical
experiments suggest that the above theoretical conclusions may translate to
practice and remain, for OFFO methods, relatively independent of noise.

Although discussed here in the context of unconstrained optimization,
adaptation of the above OFFO algorithms to problems involving convex
constraints (such as bounds on the variables) is relatively straightforward
and practical: one then needs to intersect the trust-region with the feasible
set and minimize the quadratic model in this intersection (see
\cite[Chapter~12]{ConnGoulToin00}).
It will be also of interest to further analyze the possible links between
our proposals and those of \cite{GrapStel22}, both from the theoretical and
practical perspectives,as well as to extend our investigation to the
class of adaptive regularization methods (see \cite{GratJeraToin23c}
for instance). 

{\footnotesize

\section*{\footnotesize Acknowledgements}

The authors wish to thank the referees, who suggested pointers to
additional literature and helped to significantly improve focus and
presentation. 
 
\bibliographystyle{plain}

}

\end{document}